\documentclass[11pt]{article}
\usepackage{graphicx, color}
\usepackage{amssymb}
\usepackage{amsmath}
\def\ds{\displaystyle}
\def\forall{\hbox{for all}~}
\def\L{{\bf L}}

\def\S{{\cal S}}
\def\I{{\cal I}}

\def\argmax{\hbox{arg}\!\max}

\def\bfp{{\bf p}}
\def\bft{{\bf t}}
\def\bfn{{\bf n}}

\def\ve{\varepsilon}
\def\ep{\epsilon}

\def\U{{\cal U}}

\def\R{I\!\!R}

\def\vp{\varphi}

\def\vs{\vskip 2em}

\def\v{\vskip 1em}
\def\O{{\cal O}}

\def\begi{\begin{itemize}}
\def\endi{\end{itemize}}
\def\F{{\cal F}}
\def\K{{\cal K}}
\def\C{{\cal C}}

\def\ov{\overline}
\def\Tilde{\widetilde}
\def\Hat{\widehat}
\def\wto{\rightharpoonup}
\def\bega{\begin{array}}
\def\enda{\end{array}}
\def\meas{\hbox{meas}}

\def\bel{\begin{equation}\label}
\def\eeq{\end{equation}}
\def\sqr#1#2{\vbox{\hrule height .#2pt
\hbox{\vrule width .#2pt height #1pt \kern #1pt
\vrule width .#2pt}\hrule height .#2pt }}
\def\square{\sqr74}
\def\endproof{\hphantom{MM}\hfill\llap{$\square$}\goodbreak}

\newtheorem{theorem}{Theorem}[section]
\newtheorem{lemma}[theorem]{Lemma}

\newtheorem{definition}[theorem]{Definition}

\newtheorem{example}[theorem]{Example}

\newcommand{\diff}[2]{\frac{d #1}{d#2}}
\newcommand{\pdiff}[2]{\frac{\partial #1}{\partial#2}}

\begin{document}

\title{\bf Competition Models for Plant Stems}
\vs

\author{Alberto Bressan$^{(*)}$, Sondre T.~Galtung$^{(**)}$, Audun Reigstad$^{(**)}$, and Johanna Ridder$^{(*)}$\\
\,
\\
$^{(*)}$ Department of Mathematics, Penn State University \\
University Park, Pa.~16802, USA.\\
$^{(**)}$ Department of Mathematical Sciences, \\
NTNU -- Norwegian University of Science and Technology \\ NO-7491 Trondheim, Norway.
\\
\,
\\
e-mails: axb62@psu.edu, sondre.galtung@ntnu.no, \\
audun.reigstad@ntnu.no,  johanna@jomichri.de.
}
\maketitle
\begin{abstract} The models introduced in this paper describe a uniform distribution of 
plant stems competing for sunlight. The shape of each stem, and the density of leaves, 
are designed in order to maximize the captured sunlight, subject to a  cost for transporting 
water and nutrients from the root to all the leaves.  
Given the intensity of light, 
depending on the height above ground, we first solve the optimization problem 
determining the best possible shape for a single stem.  We then study a competitive
equilibrium among a large number of similar plants, where the shape of each stem
is optimal given the shade produced by all others.  Uniqueness of equilibria is proved
by analyzing the two-point boundary value problem for a system of ODEs derived 
from the necessary conditions for optimality.
\end{abstract}

         

\vs
\section{Introduction}
\label{s:1}
\setcounter{equation}{0}
Optimization problems for tree branches have recently been studied in \cite{BPSu, BS}. In these models,
optimal shapes maximize the total amount of sunlight gathered by the leaves, subject to a cost for
building a network of branches that will bring water and nutrients from the root to all the leaves. 
Following \cite{BCM, G, MMS, X3, X15}, this cost is
defined in terms of a ramified transport.

In the present paper we consider a competition model, where a large number of similar plants compete for sunlight.
To make the problem tractable, instead of a tree-like structure we assume that each plant consists of a single stem.
As a first step, assuming that the intensity of light $I(\cdot)$ depends only on the height above ground, we determine
the corresponding optimal shape of the stem.   This will be a curve $\gamma(\cdot)$ which can be
found by classical techniques of the Calculus of Variations or optimal control \cite{BP, Cesari, FR}.
In turn, given the  density of plants (i.e., the average number of plants growing per unit area), 
if all stems have the same shape $\gamma(\cdot)$
one can compute the intensity of light $I(h)$ that reaches a point at height $h$.

An equilibrium configuration is now defined as a fixed point of the composition of the two maps 
$I(\cdot)\mapsto \gamma(\cdot)$ and $\gamma(\cdot) \mapsto I(\cdot)$.   
A major goal of this paper is to study the existence and properties of these equilibria, where the shape of each stem is
optimal subject to the presence of all other competing plants.
%The remainder of the paper is organized as follows.

In Section~2 we introduce our two basic models.    In the first model, the length $\ell$ of the stems and the 
thickness  (i.e., the density of leaves
along each stem) are assigned a priori.  The only function to optimize is thus 
the curve $\gamma:[0,\ell]\mapsto \R^2$ 
describing the shape of the stems.    
In the second model, also the length and the thickness of the stems are allowed to vary, 
and optimal values for these variables need to be determined.

In Section~3, given a light intensity function $I(\cdot)$,  we study the optimization problem for Model~1,
proving  the existence of an optimal solution and deriving necessary conditions for optimality.
We also give a condition which guarantees the uniqueness of the optimal solution. A counterexample shows that,
in general, if this condition is not satisfied multiple solutions can exist.   
In Section~4 we consider the competition of a large number of stems,
and prove the existence of an equilibrium solution. In this case, the common shape of the plant stems 
can be explicitly determined by solving a particular ODE.

The subsequent sections extend the analysis to a more general setting (Model~2), where both the length and 
the thickness of the stems are to be optimized.     In Section~5 we prove the existence of optimal stem configurations,
and derive necessary conditions for optimality, while in
Section~6 we establish  the existence of a unique equilibrium solution for the competitive game, assuming that the density 
(i.e., the average number of stems growing per unit area) is sufficiently small.
The key step in the proof is the analysis of
a two-point boundary value problem, for a system of ODEs derived 
from the necessary conditions.

In the above models, the density of stems was assumed to be uniform on the whole space.
As a consequence, the light intensity $I(h)$ depends only of the height $h$ above ground. 
Section~7, on the other hand, 
is concerned with  a family of stems growing only on the positive half line.
In this case the light intensity $I=I(h,x)$ depends also on the spatial location $x$, 
and the analysis becomes considerably more difficult. Here we only
derive a set of equations describing the competitive equilibrium,
and sketch what we conjecture should be the  corresponding shape of stems.  

The final section contains some concluding remarks.  In particular,
we discuss the issue of phototropism, i.e.~the tendency
of plant stems to bend in the direction of the light source.
Devising a mathematical model, which 
demonstrates phototropism as an advantageous trait, remains a challenging
open problem.
For a biological perspective on plant growth we refer to~\cite{LD}. 
A recent mathematical study of the stabilization problem for growing stems can be found in \cite{ABGS}.

\section{Optimization problems for a single stem}
\label{s:15}
\begin{figure}[ht]
\centerline{\hbox{\includegraphics[width=8cm]{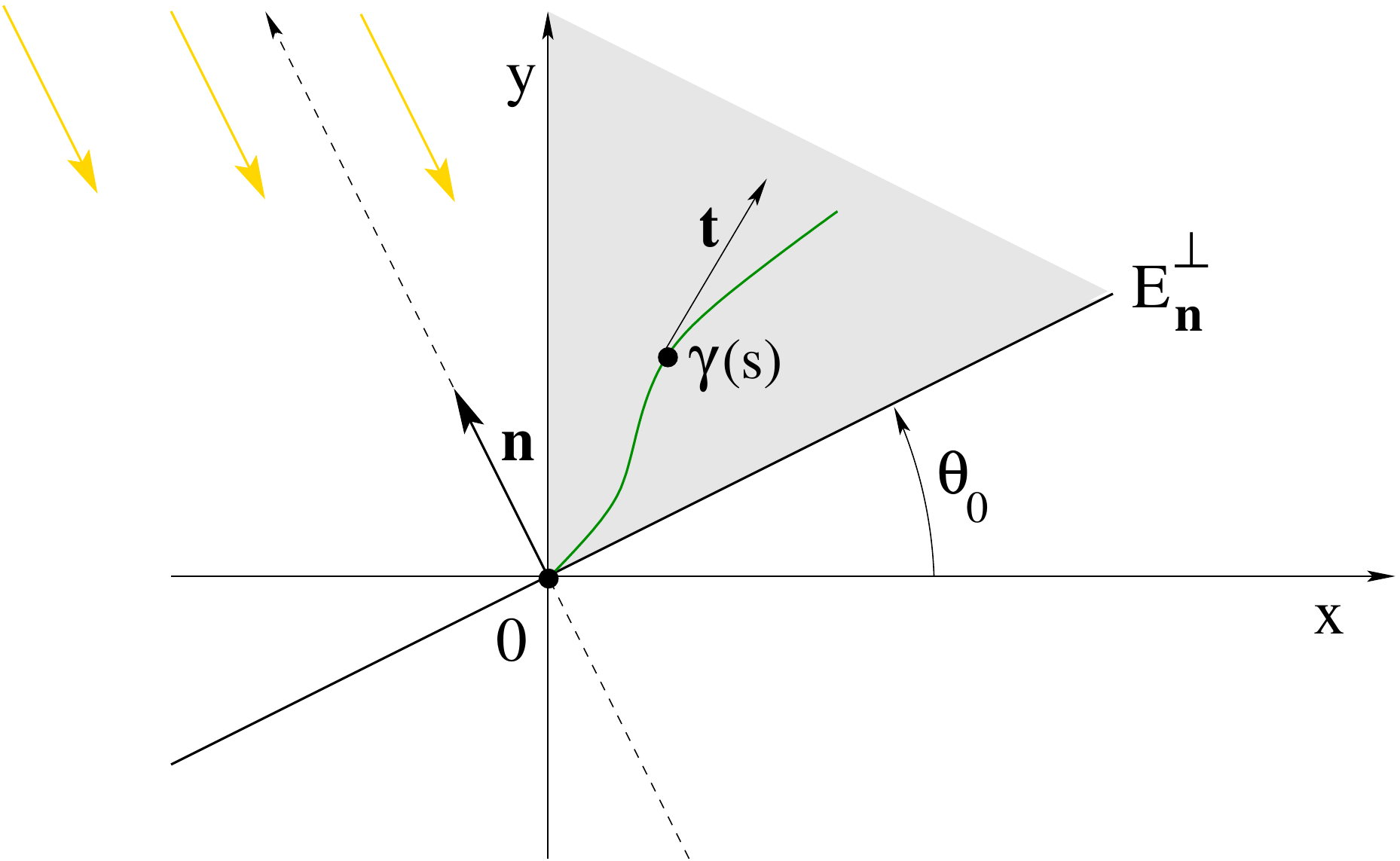}}}
\caption{\small  By a reflection argument, it is not restrictive to assume that the tangent  vector
$\bft(s)$ to the stem satisfies (\ref{thh}), i.e., it lies in the shaded cone.}
\label{f:pg43}
\end{figure}

\setcounter{equation}{0}
We shall consider plant stems
in the $x$-$y$ plane, where $y$ is the vertical coordinate.
We assume that sunlight comes from the direction of the unit vector
$$\bfn~=~(n_1, n_2), \qquad n_2<0<n_1.$$  
As in Fig.~\ref{f:pg43}, we denote by $\theta_0\in \,]0, \pi/2[\,$ the angle such that
\bel{th0} (-n_2, n_1)~=~(\cos\theta_0,\,\sin\theta_0).\eeq
Moreover, we assume that
the light intensity $I(y)\in [0,1]$ is a non-decreasing function of the height $y$.
This is due to the presence of competing vegetation: close to the ground, less light can get through.
\v
{\bf Model 1 (a stem with fixed length and constant thickness).}  
We begin by studying a simple model, where each stem has a fixed
 length $\ell$.
 Let $s\mapsto \gamma(s) = (x(s), y(s))$, $s\in [0,\ell]$, be an  arc-length parameterization
of the stem. 
As a first approximation, we assume that the leaves are uniformly distributed along the stem,
with density $\kappa$.  
The total distribution of leaves in space is thus described by a measure $\mu$, with
\bel{m1} 
\mu(A)~=~\kappa\cdot\meas\Big( \bigl\{s\in [0,\ell]\,;~~\gamma(s)\in A\bigr\}\Big)\eeq
for every Borel set $A\subseteq \R^2$.

Among all stems with given   length $\ell$,
we seek the shape which will collect the most sunlight.  This can be formulated as an optimal control problem.  Indeed,
by the Lipschitz continuity of $\gamma(\cdot)$,  the tangent vector
$$\bft(s)~=~\dot \gamma(s)~=~(\cos\theta(s), \sin\theta(s))$$
is well defined for a.e.~$s\in [0,\ell]$.  The map $s\mapsto \theta(s)$ 
will be regarded as a control function.

According to the model in \cite{BS},  calling $\Phi(\cdot)$ the density 
of the projection of $\mu$ on the  space $E^\perp_\bfn$ orthogonal to $\bfn$,
the total sunlight captured by the stem is
\bel{Sg}
\S(\gamma)~=~\int \Big( 1-\exp\{ -\Phi(z)\}\Big)\, dz~=~\int_0^\ell I(y(s))\cdot\left( 1-\exp\Big\{
{-\kappa\over\cos(\theta(s)-\theta_0)}\Big\}\right) \cos(\theta(s)-\theta_0)\, ds.
\eeq
In order to maximize (\ref{Sg}), we claim that
it is not restrictive to assume that the angle satisfies
\bel{thh}
\theta_0~\leq~\theta(s)~\leq~{\pi\over 2}\qquad\quad\forall s\in [0,\ell].\eeq
Indeed, for any measurable map $s\mapsto \theta(s)\in \,]-\pi, \pi]$, we can define a modified angle function
$\theta^\sharp (\cdot)$ by setting
\bel{ath}
\theta^\sharp (s)~=~\left\{\bega{rl}\theta(s)\quad\hbox{if} \quad \theta(s)\in \,]0, \theta_0 +\pi/2],\\[3mm]
-\theta(s)\quad\hbox{if} \quad \theta(s)\in \,]-\pi,\,\theta_0-\pi/2],\\[3mm]
2\theta_0+\pi-\theta(s)\quad\hbox{if} \quad \theta(s)\in \,]\theta_0+\pi/2, \pi],\\[3mm]
2\theta_0 -\theta(s)\quad\hbox{if} \quad \theta(s)\in \,]\theta_0-\pi/2, 0].
\enda\right.\eeq
Calling $\gamma^\sharp :[0,\ell]\mapsto \R^2$ the curve whose tangent vector is 
$\dot \gamma^\sharp (s)=(\cos\theta^\sharp (s), \sin\theta^\sharp  (s))$, since the light intensity function
$y\mapsto I(y)$ is nondecreasing we have
$\S(\gamma^\sharp )\geq\S(\gamma)$.  

By this first step, without loss of generality we can now assume $\theta(s)\in \,]0,\,\theta_0 + \pi/2]$.
To proceed further, consider the piecewise affine map
\bel{ats}
\vp(\theta)~=~\left\{\bega{cl}\theta\quad &\hbox{if} \quad \theta\in \,]\theta_0, \pi/2],\\[3mm]
\pi-\theta\quad&\hbox{if} \quad \theta\in \,[\pi/2,\, \theta_0+\pi/2],\\[3mm]
2\theta_0 -\theta\quad&\hbox{if} \quad \theta\in \,[0, \theta_0].
\enda\right.\eeq
Call $\gamma^\vp$ the curve whose tangent vector is $\dot \gamma^\vp(s)=\Big(\cos(\vp (\theta(s))), 
\,\sin(\vp(\theta (s)))\Big)$.
Since $I(\cdot)$ is nondecreasing, we again have
$\S(\gamma^\vp)~\geq~\S(\gamma)$.  
   We now observe that, since  $0<\theta_0<\pi/2$, there exists
an integer $m\geq 1$ such that the $m$-fold composition $\vp^m\doteq \vp\circ\cdots\circ \vp$
maps $[0, \theta_0+\pi/2]$ into $[\theta_0, \pi/2]$. An inductive argument 
now yields $\S(\gamma^{\vp^m})~\geq~\S(\gamma)$, completing the proof of our claim.
\v

As shown in Fig.~\ref{f:pg40}, left,  we call $ z$  the coordinate along the space $E^\perp_\bfn$ 
perpendicular to $\bfn$,
and let $y$ be the vertical coordinate. Hence
\bel{dzy}
dz(s)~=~\cos(\theta(s)-\theta_0)\, ds,\qquad\qquad dy(s)~=~\sin(\theta(s))\, ds.\eeq
In view of (\ref{thh}), one can express both $\gamma$ and $\theta$  as functions of the variable $y$. 
 Introducing the function
\bel{g}
g(\theta)~\doteq~\left(
 1-\exp\Big\{{-\kappa\over \cos(\theta-\theta_0)}\Big\}\right){\cos(\theta-\theta_0)
 \over\sin\theta}\,,\eeq
the problem can be equivalently formulated as follows.  
\begi
\item[{\bf (OP1)}] {\it Given a length $\ell>0$, find $h>0$ and a control function 
$y\mapsto\theta(y)\in [\theta_0, \,\pi/2]$ which maximizes the integral
\bel{Imax} \int_0^h I(y)\, g(\theta(y))\, dy\eeq
subject to}
\bel{ct1} \int_0^h {1\over\sin\theta(y)}\, dy~=~\ell.\eeq
\endi

\v
{\bf Model 2 (stems with variable length and thickness).}   Here we still assume that the plant consists of a single stem,
parameterized by arc-length:  $s\mapsto\gamma(s)$, $s\in [0,\ell]$.  However, now we
give no constraint on the length $\ell$ of the stem, and 
we allow the density of leaves to be variable along the stem.  

Call $u(s)$ the density of leaves at the point $\gamma(s)$.  
In other words, $\mu$ is now the measure which is absolutely continuous
w.r.t.~arc-length measure on $\gamma$, with density $u$.
Instead of (\ref{m1}) we thus have
\bel{m2} \mu(A)~=~\int_{\gamma(s)\in A} u(s)\, ds\,.\eeq

Calling $I(y)$ the intensity of light at height $y$,
the total sunlight gathered by the stem
is now computed by
\bel{S8} \S(\mu)~=~\int_0^\ell I(y(s))\cdot\left( 1-\exp\Big\{
{-u(s)\over\cos(\theta(s)-\theta_0)}\Big\}\right) \cos(\theta(s)-\theta_0)\, ds
\eeq
As in \cite{BS}, we consider a cost for transporting water and nutrients from the root to the leaves.
This is measured by
\bel{TC}
\I^\alpha(\mu)~=~\int_0^\ell \left(\int_s^\ell u(t)\, dt\right)^\alpha ds,\eeq
for some $0<\alpha<1$.   Notice that, in Model 1, this cost was 
the same for all stems and hence it
did not play a role in the optimization.

For a given constant $c>0$, we now consider a second optimization problem:
\bel{OP2}
\hbox{maximize:}\quad \S(\mu)- c \I^\alpha (\mu),\eeq
subject to:
\bel{C6} y(0)~=~0,\qquad\qquad \dot y(s)~=~\sin \theta(s).\eeq
The maximum is sought over all controls $\theta:\R_+\mapsto [0, \pi]$
and $u:\R_+\mapsto \R_+$.
Calling
\bel{zd3}z(t)~\doteq~\int_t^{+\infty} u(s)\, ds,\eeq
\bel{gtu}G(\theta,u)~\doteq~\left( 1 - \exp\left\{{-u\over \cos(\theta-\theta_0)}
\right\}\right)\cos(\theta-\theta_0)\,,\eeq
this leads to an optimal control problem in a more standard form.
\begi
\item[{\bf (OP2)}] 
{\it Given a sunlight intensity function $I(y)$, and constants $0<\alpha<1$,
$c>0$, find controls $\theta:\R_+\mapsto [\theta_0, \pi/2]$
and $u:\R_+\mapsto \R_+$ which maximize the integral
\bel{max33} \int_0^{+\infty} \Big[ I(y)\, G(\theta, u) - c\,z^\alpha\Big]\, dt ,\eeq
subject to
\bel{max44} 
\left\{\bega{rl} \dot y(t)&=~\sin\theta ,\\[3mm]
\dot z(t)&=~- u,\enda\right.
\qquad\qquad
\left\{\bega{rl} y(0)&=~0,\\[3mm]
z(+\infty)&=~0.\enda\right.
\eeq
}
\endi

\section{Optimal stems with fixed length and thickness}
\label{s:16}
\setcounter{equation}{0}

\subsection{Existence of an optimal solution.}

Let  $I(y)$  be the light intensity, which we assume is a non-decreasing 
function of the vertical component $y$.
For a given $\kappa>0$ (the thickness of the stem), we seek a curve $s\mapsto \gamma(s)$,
starting at the origin and with a fixed length $\ell$,  which maximizes the sunlight
functional defined at (\ref{S8}).  \v

\begin{theorem}\label{t:31}
 For any non-decreasing function $y\mapsto I(y)\in [0,1]$ and any 
constants $\ell,\kappa>0$ and $\theta_0\in\,]0,\pi/2[\,$, the optimization problem
{\bf(OP1)} has at least one solution. 
\end{theorem}
\v
{\bf Proof.} {\bf 1.} Let $M$ be the supremum among all admissible payoffs in (\ref{Imax}). 
By the analysis in \cite{BS} it follows that
$$0~\leq~M~\leq~\kappa\,\mu(\R^2)~=~\kappa\,\ell.$$
Hence there exists a maximizing sequence of control functions 
$\theta_n:[0, h_n]\mapsto [\theta_0,\pi/2]$, so that
\bel{ctn} \int_0^{h_n} {1\over\sin\theta_n(y)}\, dy~=~\ell\qquad\forall n\geq 1,\eeq
\bel{JN} \int_0^{h_n} I(y) g(\theta_n(y))\, dy~\to~M.\eeq
\v
{\bf 2.} For each $n$, let $\theta^\sharp_n$ be the non-increasing rearrangement of 
the function $\theta_n$.  
Namely, $\theta^\sharp_n$ is the unique
(up to a set of zero measure) non-increasing function such that, for every $c\in \R$
\bel{eme}\meas\Big(\{s\,;~~\theta^\sharp_n(s)< c\}\Big)~=~\meas\Big(\{s\,;~~\theta_n(s)< c\}\Big).\eeq
This can be explicitly defined as
$$\theta^\sharp_n(y)~=~\sup \left\{  \xi\,;~~ \meas\Big(\{\sigma\in [0,h_n]\,;~~
\theta_n(\sigma)
\geq \xi\}\Big)>y\}\right\}.$$
For every $n\geq 1$ we claim that
\bel{rear1} \int_0^{h_n} {1\over\sin\theta^\sharp_n(y)}\, dy~=~\int_0^{h_n} {1\over\sin\theta_n(y)}\, dy~=~\ell,\eeq
\bel{rear2}\int_0^{h_n} I(y) g(\theta^\sharp_n(y))\, dy~\geq~
\int_0^{h_n} I(y) g(\theta_n(y))\, dy.\eeq
Indeed, to prove the first identity we observe that, by (\ref{eme}), there exists
a measure-preserving map $y\mapsto \zeta(y)$ from $[0, h_n]$ into itself such that
$\theta_n^\sharp(y)~=~\theta_n(\zeta(y))$.   Using $\zeta$ as new variable of integration,
one immediately obtains (\ref{rear1}).

To prove (\ref{rear2}) we observe that the function $g$ introduced at (\ref{g}) is smooth 
and satisfies 
\bel{g'} g'(\theta)~\leq~0\qquad\forall \theta\in [\theta_0, \,\pi/2].\eeq
Therefore, the map $y\mapsto g(\theta_n^\sharp(y))$ coincides with the non-decreasing 
rearrangement of $y\mapsto g(\theta_n(y))$.
On the other hand, since $I(\cdot)$ is non-decreasing, it trivially 
coincides with the non-decreasing 
rearrangement of itself.   Therefore,  (\ref{rear2}) is an immediate consequence 
of  the Hardy-Littlewood inequality
\cite{LL}. 
\v
{\bf 3.} Since all functions $\theta_n^\sharp$ are non-increasing, they have 
bounded variation.   Using Helly's compactness theorem, by possibly  extracting a subsequence, we can find $h>0$ and 
a non-increasing function $ \theta^*:[0, h]\mapsto [\theta_0, \pi/2]$ such that
\bel{conv} \lim_{n\to\infty} h_n~=~ h\,,
\qquad\qquad  \lim_{n\to\infty} \theta_n^\sharp(y) ~=~ \theta^*(y)
\qquad\hbox{for a.e.}~y\in [0,h].\eeq
This implies
$$\ \int_0^{h} {1\over\sin \theta^*(y)}\, dy~=~\ell,\qquad\qquad
\int_0^{h} I(y) g( \theta^*(y))\, dy~=~M,$$
proving the optimality of $ \theta^*$.
\endproof

%{\bf Remark.} For $\ve>0$ small, a leading order approximation is obtained by 
%replacing the function $1-e^{-\Phi(z)}$ with $\Phi(z)-\Phi^2(z)$.
%This may lead to simpler necessary conditions for optimality.

\subsection{Necessary conditions for optimality}
Let  $y\mapsto \theta^*(y) $ be an optimal solution.   By the previous analysis we already 
know that the function $\theta^*(\cdot)$ is non-increasing.   Otherwise, its non-increasing rearrangement achieves a better payoff.
In particular, this implies that the left limit at the terminal point $y=h$ is well defined:
\bel{tll}\theta^*(h)~=~\lim_{y\to h-} \theta^*(y).\eeq

Consider an arbitrary perturbation 
$$\theta_\epsilon ~=~\theta^* + \epsilon \Theta, \qquad\qquad h_\ep~=~h+\ep\eta.$$
The constraint (\ref{ct1}) implies
\bel{nct} \int_0^{h+\ep\eta} {1\over \sin\theta_\ep(y)}\, dy~=~\ell.\eeq
Differentiating (\ref{nct}) w.r.t.~$\ep$
one obtains
\bel{nc5} {1\over \sin\theta^*(h)}\,\eta - \int_0^h {\cos\theta^*(y)\over \sin^2\theta^*(y)}\, 
\Theta(y)\,dy
~=~0.\eeq
Next, calling 
$$J_\epsilon~\doteq~\int_0^{h_\epsilon} I(y)g(\theta_\epsilon(y)) dy$$
and assuming that $I(\cdot)$ is continuous at least at $y=h$, by (\ref{nc5}) we obtain
\bel{NC4} \bega{l}
0~=~\ds {d\over d\epsilon} J_\epsilon\bigg|_{\epsilon=0}~=~\ds \int_0^h I(y)
%\left( {g'(\theta^*(y))\over \sin \theta^*(y)}
% - {g(\theta^*(y))\cos \theta^*(y)\over \sin^2 \theta^*(y) } \right)
g'(\theta^*(y))\Theta(y)\, dy
 \\[4mm]
 \qquad\qquad\qquad \qquad\qquad\ds
 + I(h) g(\theta^*(h))\cdot \sin \theta^*(h) \int_0^h {\cos
 \theta^*(y )\over \sin^2 \theta^*(y)} \,\Theta(y)\, dy.
 \enda\eeq
Since (\ref{NC4}) holds for arbitrary perturbations $\Theta(\cdot)$, 
the optimal control 
$\theta^*(\cdot)$ 
should satisfy the identity
\bel{NC0}
I(y) g'\bigl(\theta^*(y)\bigr)
+\lambda \cdot{\cos\theta^*(y)\over \sin^2 \theta^*(y)}~=~0,
\qquad\quad\hbox{for a.e.}~y\in [0,h], \eeq
where
\bel{lam}
\lambda~=~ I(h) g(\theta^*(h))\cdot \sin \theta^*(h).\eeq
It will be convenient to write
\bel{gG}g(\theta)~=~{G(\theta)\over\sin\theta}\,,\qquad\quad G(\theta)~\doteq~
\left(
 1-\exp\Big\{{-\kappa\over \cos(\theta-\theta_0)}\Big\}\right)\cos(\theta-\theta_0).
\eeq
Inserting (\ref{gG}) in (\ref{NC0}) one obtains the pointwise identities
\bel{NC5}
I(y)\Big( {G'(\theta^*(y))\sin \theta^*(y)}
 - {G(\theta^*(y))\cos\theta^*(y)} \Big)+\lambda\cdot
{\cos\theta^*(y)}~=~0,\eeq
%$$
% G'(\theta^*(y))\,\tan \theta^*(y)
% - G(\theta^*(y)) ~= ~{\lambda\over I(y)}\,.$$
At $y=h$, the identities (\ref{lam}) and (\ref{NC5})  yield
$$ G'(\theta^*(h))\,\tan \theta^*(h)
 - G(\theta^*(h)) ~= ~{ I(h) G(\theta^*(h))\over I(y)}\,.$$
Hence
$$G'(\theta^*(h))\,\tan \theta^*(h)~=~0,$$
which implies 
\bel{th*}\theta^*(h)~=~\theta_0\,,\qquad\qquad \lambda~=~I(h) g(\theta_0) \sin\theta_0
~=~\bigl(
 1-e^{-\kappa}\bigr)\,I(h)\,.\eeq

Notice that (\ref{NC5}) corresponds to
\bel{NC6}
\theta^*(y)~=~\argmax_{\theta\in[0,\pi]}~\left\{ I(y) {G(\theta)\over \sin \theta}
+{\lambda\over \sin \theta}\right\}\,.\eeq
Equivalently, $\theta=\theta^*(y)$ is the solution to
\bel{NC7}
G'(\theta)\tan\theta - G(\theta)~=~-{\lambda\over I(y)}\,,\eeq
where $G$ is the function at (\ref{gG}).

\begin{lemma}\label{l:vp}
Let $G$ be the function at (\ref{gG}).   Then for every $z\in \,]-\infty, \,e^{-\kappa}-1]$
the equation
\bel{ET} F(\theta)~\doteq~G'(\theta)\tan\theta - G(\theta)~=~z\eeq
has a unique solution  $\theta=\vp(z)\in [\theta_0, \pi/2[\,$.
\end{lemma}
\v
%{\color{red} What else can we say about this function $\vp$ ?    This will be useful in order to understand uniqueness of the equilibrium solution.   See the Cauchy problem (\ref{BCP}).}\v
{\bf Proof.}
Observing that 
\bel{Gprop}
\left\{ \bega{rl} G(\theta_0)&=~1-e^{-\kappa},\\[3mm]
 G'(\theta_0)&=~0,\enda\right.
 \qquad\qquad \left\{\bega{rl} G'(\theta)&<~0\\[3mm] 
 G''(\theta)&<~0\enda\right.\qquad 
 \hbox{for}~~\theta\in \,]\theta_0, \pi/2[\,,\eeq
  we obtain
$F(\theta_0)=e^{-\kappa}-1$ and
$$ F'(\theta)\,=\,G''(\theta)\tan\theta+G'(\theta)\tan^2\theta~<~0\qquad \hbox{for}
~~ \theta \in [\theta_0, \pi/2[\,. $$
Therefore, for $\theta\in [\theta_0, \pi/2[\,$, the left hand side of (\ref{ET})   is monotonically decreasing from $e^{-\kappa}-1$ to $-\infty$. 
We conclude  that~(\ref{ET}) has a unique solution 
$\theta=\varphi(z)$ for any $z\in \,]-\infty, \,e^{-\kappa}-1]$.  
\endproof

 The optimal control $\theta^*(\cdot)$
determined by the necessary condition (\ref{NC7}) is thus recovered by
\bel{opthe} \theta^*(y)~=~\varphi\left(\frac{-\lambda}{I(y)}\right)~=~
\varphi\left(\frac{(e^{-\kappa}-1 )I(h)}{I(y)}\right).\eeq
Next, we need to determine  $h$ so that  the constraint
\bel{Leq}
L(h)~\doteq~\int_0^h\frac{1}{\sin(\theta^*(y))}\,dy~=~\ell
\eeq
is satisfied.
As shown by Example~\ref{e:1} below, 
the solution of~(\ref{opthe})-(\ref{Leq}) may not be unique. 

In the following, we seek a condition on $I$ which implies  that $L$ is monotone, i.e., 
\bel{Lmon}
L'(h)~=~\frac{1}{\sin(\theta_0)}+\int_0^h\frac{\cos\theta^*(y)}{\sin^2\theta^*(y)}\frac{1}{F'(\theta^*(y))}\frac{I'(h)}{I(y)}G(\theta_0)\,dy ~>~0\,.
\eeq
This will guarantee that (\ref{Leq}) has a unique solution. 
To get an upper bound for $F'(\theta)$, observe that, for $\theta\in [\theta_0, \pi/2[$, 
\[\bega{l} 
F'(\theta)~\le ~\ds \tan(\theta)G''(\theta)\\[3mm]
=\ds~ -\tan(\theta)\bigg[ \cos(\theta-\theta_0)\left(1-\left(\frac{\kappa}{\cos(\theta-\theta_0)}
+1\right)\exp\Big\{{-\kappa\over \cos(\theta-\theta_0)}\Big\}
\right) \\[3mm]
\qquad\qquad \qquad\qquad \qquad\qquad \qquad\qquad 
\ds+ \frac{\tan^2(\theta-\theta_0)}{\cos(\theta-\theta_0)}\kappa^2 \exp\Big\{{-\kappa\over \cos(\theta-\theta_0)}\Big\}\bigg]\\[3mm]
=~-\tan(\theta) \,\cos(\pi/2-\theta_0) \bigl( 1-(\kappa+1)e^{-\kappa}\bigr).
\enda\]
Since $\theta^*(y)\in [\theta_0, \,\pi/2]$ and $G(\theta_0) = 1-e^{-\kappa}$, 
using the above inequality one obtains
$$
 \int_0^h \frac{\cos\theta^*(y)}{\sin^2\theta^*(y)}\cdot  \frac{1}{| F'(\theta^*(y))|}\frac{I'(h)}{I(y)}G(\theta_0)\,dy~\le~ \frac{\cos^2\theta_0}{\sin^3 \theta_0}\cdot
 \frac{1-e^{-\kappa}}{\cos(\pi/2-\theta_0) \Big( 1-(\kappa+1)e^{-\kappa}\Big)}\int_0^h \frac{I'(h)}{I(y)}\,dy\,.
$$
Hence (\ref{Lmon}) is satisfied provided that 
\bel{uni}
\int_0^h \frac{I'(h)}{I(y)}\,dy~<~
 \tan^2\theta_0\cdot \frac{\cos(\pi/2-\theta_0) \bigl( 1-(\kappa+1)e^{-\kappa}\bigr)}{1-e^{-\kappa}}\,.
\eeq
{}From the above analysis, we conclude
\v
\begin{theorem}\label{t:1} Assume that the light intensity function $I$ is Lipschitz continuous and satisfies 
the strict inequality (\ref{uni}) for a.e.~$h\in [0,\ell]$.   
Then  the optimization problem {\bf (OP1)} has a unique optimal solution
$\theta^*:[0, h^*]\mapsto [\theta_0, \,\pi/2]$.
The function $\theta^*$ is non-increasing, and satisfies 
\bel{t*}
\theta^*(y)~=~\vp\left( (e^{-\kappa}-1) {I(h^*)\over I(y)}\right),\eeq
where $z\mapsto \vp(z)=\theta $ is the function implicitly defined by (\ref{ET}).
\end{theorem}
\v

\begin{figure}[ht]
\centerline{\hbox{\includegraphics[width=15cm]{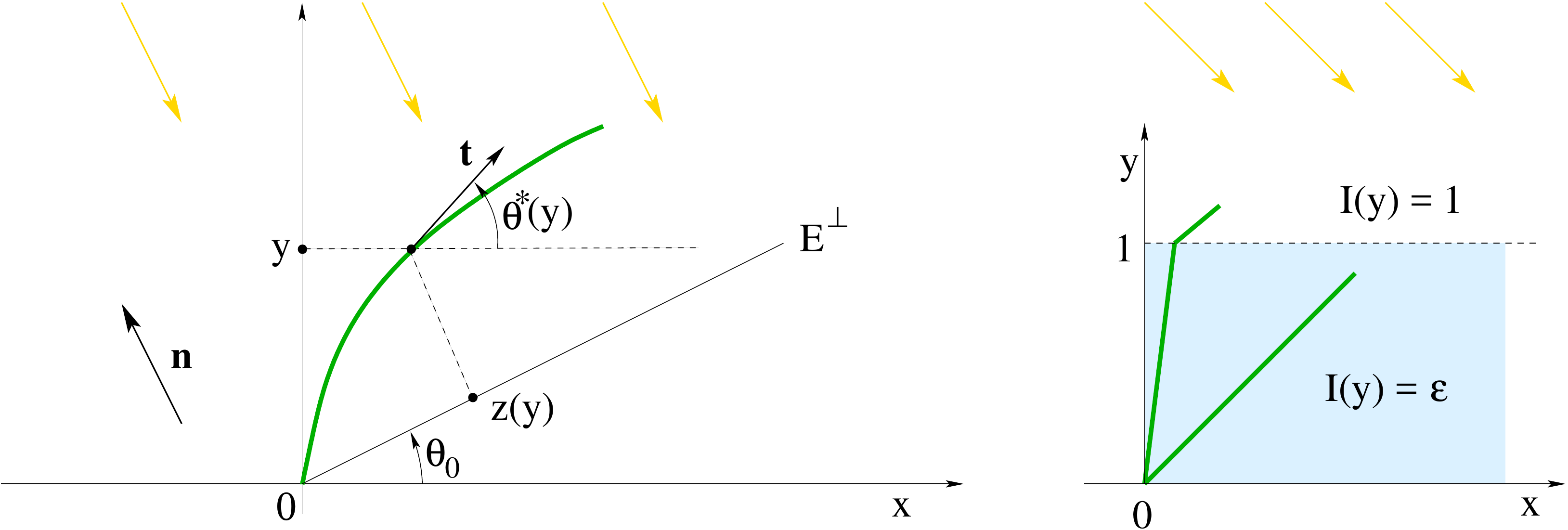}}}
\caption{\small   Left: the optimal shape of a stem, %for a given light intensity function $I(\cdot)$, 
as described in Theorem~\ref{t:1}.
Right: if the light intensity $I$ changes abruptly as a function of the hight, 
the optimal shape may not be unique, as
shown in Example~\ref{e:1}.}
\label{f:pg40}
\end{figure}

The following example shows that, without the  bound (\ref{uni}) 
on  the 
sunlight intensity function $I(\cdot)$, the conclusion of Theorem~\ref{t:1} can fail.
\v
\begin{example} \label{e:1}{\bf (non-uniqueness).}   {\rm
Choose $\bfn= \Big(-{1\over\sqrt 2}, {1\over\sqrt 2}\Big)$, ~$\ell = 6/5<\sqrt 2$, $\kappa=1$,
$$I(y)~=~\left\{\bega{cl} \ve \qquad &\hbox{if}\quad y\in [0,1],\\[3mm]
1\qquad &\hbox{if}\quad y> 1,\enda\right.$$
with $\ve>0$ small.

By Theorem~\ref{t:31} at least one optimal solution exists.
By the previous analysis, any optimal solution $\theta^*:[0, h^*]\mapsto [\theta_0, \pi/2]$ 
satisfies  
the necessary conditions
(\ref{t*}). In this particular case, this implies that $\theta^*(y)$ is constant separately
for $y<1$ and for $y>1$.   As shown in Fig.~\ref{f:pg40}, right, these necessary conditions can 
have two  solutions.
\v
{\bf Solution 1.}   If    $h^*<1$, then $I(y) = \ve$ for all $y\in [0,h^*]$ and the necessary conditions (\ref{t*})
yield
\bel{SO1}\theta_1^*(y)~=~\vp(e^{-1}-1)~=~\theta_0~=~\pi/4\qquad\forall y\in [0,h^*].$$
The total sunlight collected is
$$\S_\ve(\theta^*_1) ~=~{6\over 5} (1-e^{-1})\, \ve.\eeq

{\bf Solution 2.} If $h^*>1$, then $I(h^*)=1$ and the necessary conditions (\ref{t*})
yield
$$\theta_2^*(y)~=~\vp\left((e^{-1}-1) {I(h^*)\over I(y)}\right)
~=~\left\{ \bega{cl}\ds
\vp\left((e^{-1}-1) \ve^{-1}\right)\quad &\hbox{if}\quad y\in [0,1],\\[3mm]
 \ds {\pi/ 4} \quad &\hbox{if}\quad y>1.\enda\right.$$
 Calling $\alpha=\alpha(\ve)\doteq~\vp\left((e^{-1}-1) \ve^{-1}\right)$,
the total sunlight collected in this case is 
\bel{SO2}
\S_\ve(\theta_2^*)~=~ \left( 1 - \exp\left\{- {1 \over \cos(\alpha-\pi/4)}\right\} \right)\cos(\alpha-\pi/4) \,\ve + \left( \frac{6}{5}- {1\over \sin\alpha}\right) \bigl(1 - e^{-1}\bigr).
\eeq
We claim that, for a suitable choice of $\ve\in \,]0,1[\,$, the two quantities in (\ref{SO1})
and (\ref{SO2}) become equal.   Indeed, as $\ve\to 0+$ we have
$$\alpha(\ve)~\doteq~\vp\left({e^{-1}-1\over  \ve}\right)~\to ~
{\pi\over 2}\,,$$
\bel{se1}
\S_\ve(\theta^*_1) ~\to~0,\qquad\qquad  \S_\ve(\theta_2^*)~\to~{1 - e^{-1}\over 5}\,.\eeq
On the other hand, as $\ve\to 1$ we have $\alpha(\ve)\to \pi/4$.
By continuity, there exists $\ve_1\in \,]0,1[\,$ such that 
$$\sin \alpha(\ve_1)~=~{5\over 6}\,.$$
As $\ve\to \ve_1+$, we have
\bel{se2}\S_\ve(\theta_2^*)~\to~ \left( 1 - \exp\left\{- {1 \over \cos(\alpha(\ve_1)-\pi/4)}\right\} \right)\cos(\alpha(\ve_1)-\pi/4) \,\ve_1~<~ \S_{\ve_1}(\theta_1^*).\eeq
Comparing (\ref{se1}) with (\ref{se2}), by continuity we conclude that there exists
some $\Hat \ve\in \,]0, \ve_1[\,$ such that 
$\S_{\Hat \ve}(\theta_1^*)~=~\S_{\Hat \ve}(\theta_2^*)$.   Hence for $\ve=\Hat \ve$ 
the optimization problem
has two distinct solutions.
}
\end{example}

\v

We remark that in this example the light intensity $I(y)$ is discontinuous at $y=1$.
However,  by 
a mollification one can still construct a similar example
with two optimal configurations, also for $I(\cdot)$ smooth.
Of course, in this case the derivative $I'(h)$ will be extremely large for $h\approx 1$,
so that the assumption (\ref{uni}) fails.

%\begin{figure}[ht]
%\centerline{\hbox{\includegraphics[width=10cm]{FIG/pg34.eps}}}
%\caption{\small  If $u(y)<0$ on some interval, we can replace the curve $\gamma$
%with another curve (blue), collecting more sunlight.
%}
%\label{f:pg34}
%\end{figure}

\section{A competition model}
\label{s:44}
\setcounter{equation}{0}
In the previous analysis, the light intensity function $I(\cdot)$ was a priori given.
We now consider a continuous distribution of stems, and determine 
the average sunlight $I(y)$ available at height $y$ above ground, depending
on the density of vegetation above $y$.
 
Let the constants $\ell, \kappa>0$ be given, specifying the length and thickness of each stem.  
We now introduce another constant $\rho>0$ describing the density of stems,
i.e.~how many stems grow per unit area.   
Assume that all stems have the same height and shape, described the the function 
$\theta:[0, h]\mapsto [\theta_0, \pi/2]$.  For any $y\in [ 0,h]$,
the total amount of vegetation at height $\geq y$, per unit length, 
is then measured by 
$$\rho \cdot \int_y^h {\kappa\over \sin \theta(y)}\, dy.$$
The corresponding light intensity function is defined as
\bel{I2}
I(y)~\doteq~\exp\left\{ -\rho \cdot \int_y^h {\kappa\over \sin \theta(y)}\, dy\right\}\qquad 
\hbox{for}\quad y\in [0,h],\eeq
while $I(y)=1$ for $y\geq h$.
We are interested in equilibrium configurations, where the shape of the stems is optimal
for the light intensity $I(\cdot)$. 
We recall that  $\theta_0$ is the angle of  incoming light rays, as in (\ref{th0}),
while the constants $\ell,\kappa>0$ denote the length and thickness of the stems.
\v
\begin{definition}\label{d:E1}  Given an angle $\theta_0\in \,] 0, \pi/2]$ 
and constants $\ell, \kappa,\rho>0$, we say that 
a light intensity function $I^*:\R_+\mapsto [0,1]$ and a stem shape function 
$\theta^*:[0, h^*]\mapsto 
[\theta_0, \pi/2]$
yield a {\bf competitive equilibrium} if the following holds.
\begi
\item[(i)] The stem shape function $\theta^*:[0, h^*]\mapsto 
[\theta_0, \pi/2]$ provides an optimal solution to the optimization problem 
{\bf (OP1)}, with light intensity function $I=I^*$.
\item[(ii)] For all $y\geq 0$, the  light intensity  at height $y$ satisfies
\bel{I3} 
I^*(y)~=~\exp\left\{ -\rho \cdot \int_{\min\{y,h^*\}}^{h^*} {\kappa\over \sin \theta^*(y)}\, dy\right\}.\eeq
\endi
\end{definition}
\v
If the density of vegetation is sufficiently small,  we now  show that an equilibrium configuration exists.
\v
\begin{theorem}
\label{t:2}  Let the light angle $\theta_0\in \,]0, \pi/2]$ and the stem length $\ell>0$ be given.
Then there exists a constant $c_0>0$ such that, whenever $\kappa\,\rho\leq c_0$, 
an equilibrium configuration exists.
\end{theorem}
\v
{\bf Proof.} {\bf 1.} Let $\K$ be the set of all couples $(\bar h,\theta)$, where $\ov h\in [0,\ell]$ and
$\theta:[0,\ell]\mapsto[\theta_0, \pi/2]$ is any non-increasing function.
We observe that $\K$ is a compact, convex subset of the product space $\R\times
\L^1([0, \ell])$.

For all $y\geq 0$, define the light intensity
\bel{I4} 
I(y)~=~\exp\left\{ - \int_{\min\{y,\bar h\}}^{\bar h} {\rho\, \kappa\over \sin \theta(y)}\, dy\right\}.\eeq
For $\rho>0$ small enough, we claim that  this function satisfies the assumption
of Theorem~\ref{t:1}. Indeed,  for a.e.~$h\in [0,\bar h]$ 
the left hand side of (\ref{uni}) is estimated by
$$I'(h)\, \int_0^h  {1\over I(y)}\, dy~=~{\rho \kappa\over \sin \theta(h)} \cdot 
\int_0^h \exp\left\{\rho \int_{\min\{ y, \bar h\}}^h{\kappa\over\sin\theta(y')}\, dy' \right\}\, dy,
$$
and it clearly approaches zero as $\rho\to 0$.   On the other hand, for $h>\bar h$
we have $I'(h)=0$, hence the inequality (\ref{uni}) is trivially satisfied.

By Theorem~\ref{t:1}, the optimization problem {\bf (OP1)} has  a unique solution 
$\theta^*:[0, h^*]\mapsto [\theta_0, \pi/2]$.   
For convenience, we extend this map to the entire interval $[0,\ell]$ by setting
\bel{t8}\theta^*(y)=\theta_0\qquad \qquad y\in [h^*, \ell]\,.\eeq

The above definition yields a mapping 
\bel{Lade}\Lambda: (\bar h, \theta)~\mapsto~(h^*,\theta^*)\eeq
from $\K$ into itself.     
\v
{\bf 2.} We claim that the map $\Lambda$ in (\ref{Lade}) is continuous.   

Otherwise, there would exist a sequence $(\bar h_n,\theta_n)\to (\bar h,\theta)$,
such that
$\Lambda(\bar h_n,\theta_n)$ does not converge to $\Lambda (\bar h,\bar\theta)$.
By compactness we can choose a subsequence $(\bar h_{n_k},\theta_{n_k})$
such that 
$$\Lambda(\bar h_{n_k},\theta_{n_k})~\to~(h^\sharp,\theta^\sharp)~\not=~
\Lambda(\bar h,\theta)$$
But then $(h^\sharp,\theta^\sharp)$ would be a second optimal solution to the optimization
problem with light intensity function (\ref{I4}), in contradiction with the uniqueness stated in Theorem~\ref{t:1}.

 By Schauder's theorem, the continuous map $\Lambda$
from the compact convex set $\K$ into itself has a fixed point,
which provides the desired equilibrium configuration.
\endproof

\begin{figure}[ht]
\centerline{\hbox{\includegraphics[width=10cm]{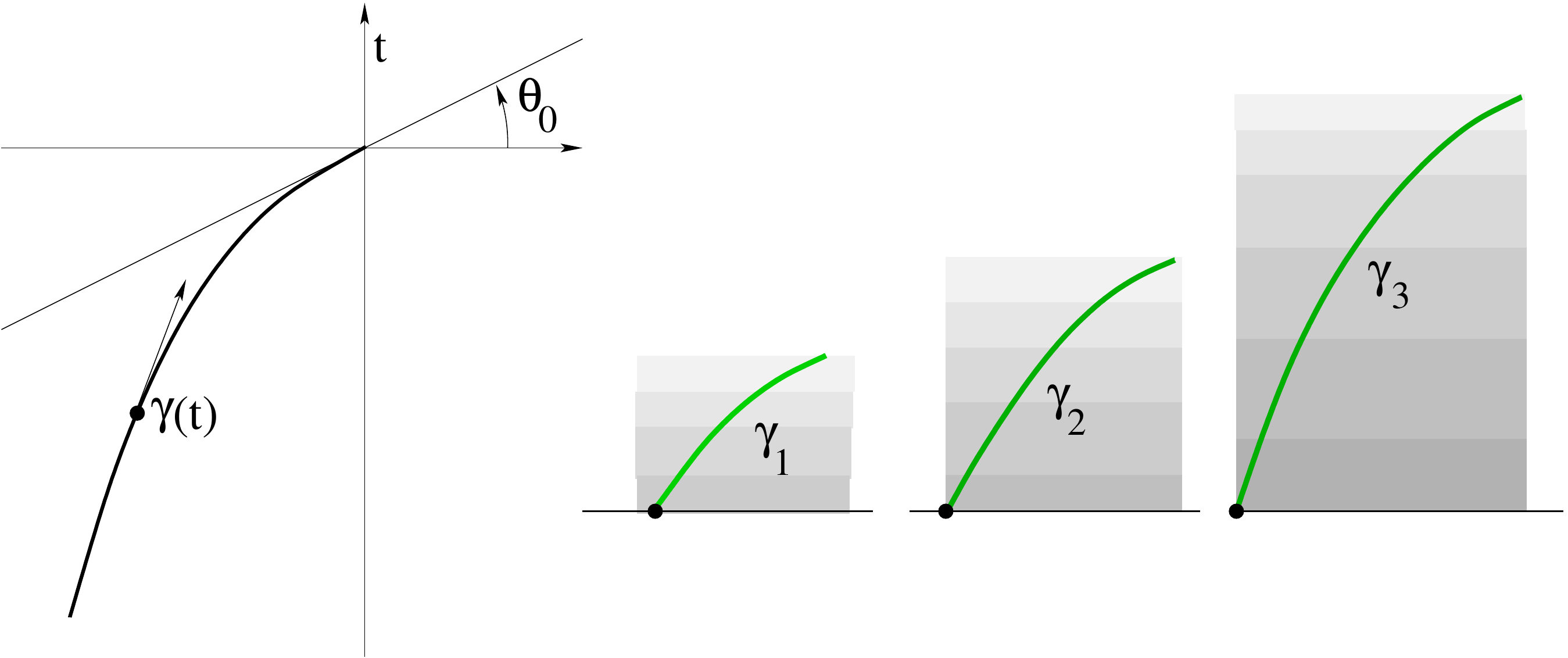}}}
\caption{\small  Left: the  curve $\gamma$, parameterized by the coordinate $t$.
For $t<0$, the tangent vector is $ {d\gamma\over dt} = (\tan \theta(t), 1)$, where $\theta(t)$
is obtained  by solving the Cauchy problem  (\ref{BCP}).
Right: for different lengths $0<\ell_1<\ell_2<\ell_3$, the equilibrium configuration is obtained by taking the upper portion of the same curve $\gamma$, up to the length $\ell_i$, $i=1,2,3$.}
\label{f:ir89}
\end{figure}

\subsection{Uniqueness and representation of equilibrium solutions.}
By (\ref{opthe}) and  (\ref{I3}),   this
equilibrium configuration $(h^*, \theta^*)$ must satisfy the necessary condition
\bel{NE1}
\theta^*(y)~=~\vp\left((e^{-\kappa}-1) \exp\left\{\int_{y}^{h^*} {\rho\kappa\over \sin \theta^*(y)}\, dy
\right\}\right),\qquad\qquad y\in [0, h^*],\eeq
where $\vp$ is the function defined in Lemma~\ref{l:vp}.
Here the constant $h^*$ must be determined so that
\bel{co4}\int_0^{h^*} {1\over \sin \theta^*(y)}\, dy~=~\ell.\eeq

Based on (\ref{NE1}), one  obtains a simple representation of  all equilibrium configurations, for any 
length $\ell>0$.  Indeed,  for $t\in \,]-\infty, 0]$,
let $t\mapsto \Hat \zeta(t)$ be the solution of the Cauchy problem
%\bel{odee}
$$\zeta'~=~-{\rho \kappa\over \sin \theta},\qquad
 \hbox{where}
\qquad \theta~=~\vp\Big( (e^{-\kappa}-1)\,e^\zeta\Big),
$$
with terminal condition
$ \zeta(0)\,=\,0$.

Notice that the corresponding 
function $t\mapsto \Hat\theta(t)= \vp\Big( (e^{-\kappa}-1)\,e^{\Hat\zeta(t)}\Big)$   satisfies
$$\Hat\theta(0)~=~\vp(e^{-\kappa}-1) ~=~\theta_0\,.$$

For any length $\ell$ of the stem, choose $h^*=h^*(\ell)$ so that
\bel{h**}\int_{-h^*}^0 {1\over \sin \Hat\theta(t)}\, dt~=~\ell\,.\eeq
The shape of the stem that achieves the competitive equilibrium is then
provided by
\bel{t**} \theta^*(y)~=~\Hat\theta(y-h^*)\,,\qquad\qquad y\in [0, h^*].\eeq

Since the backward Cauchy problem
\bel{BCP} \zeta'~=~-{\rho \kappa\over \sin\Big(\vp\bigl( (e^{-\kappa}-1)\,e^\zeta\bigr)\Big)}\,,
\qquad\qquad \zeta(0)\,=\,0,\eeq
has a unique solution,  we conclude that, if an equilibrium solution exists, by the representation (\ref{t**}) it must be unique.

\section{Stems with variable length and thickness}
\label{s:5}
\setcounter{equation}{0}

We now consider the optimization problem {\bf (OP2)}, allowing for  stems of different
lengths and with variable density of leaves.

\subsection{Existence of an optimal solution.}

\begin{theorem}\label{t:51}
 For any bounded, non-decreasing function $y\mapsto I(y)\in [0,1]$ and any 
constants $0<\alpha<1$, $c>0$ and $\theta_0\in\,]0,\pi/2[\,$, the optimization problem
{\bf(OP2)} has at least one solution. 
\end{theorem}
\v
{\bf Proof.}  {\bf 1.}  Consider a maximizing sequence of couples
$(\theta_k, u_k): \R_+\mapsto  [\theta_0, \pi/2] \times \R_+$.    For $k\geq 1$, 
let
$$s~\mapsto ~\gamma_k(s)~=~\left(\int_0^s \cos\theta_k(s)\,ds\,,~\int_0^s \sin\theta_k(s)\,ds
\right)$$
be the arc-length parameterization of the stem $\gamma_k$.
Call $\mu_k$ the Radon measure on $\R^2$ describing the distribution of leaves along 
$\gamma_k$.   For every Borel set $A\subseteq\R^n$, we thus have
\bel{muk2}
\mu_k(A)~=~\int_{\gamma_k(s)\in A}~u_k(s)\, ds.\eeq
For a given radius $\rho>0$, we have the decomposition 
$$\mu_k~=~ \mu_k^\flat+\mu_k^\sharp\,,$$ where
$\mu_k^\flat$ is the restriction of $\mu_k$ to the ball $B(0,\rho)$, while 
$\mu_k^\sharp$ the restriction of $\mu_k$ to the complement  $\R^2\setminus B(0,\rho)$.
By the same arguments used in steps {\bf 1-2}  of the proof of Theorem~3.1 in \cite{BPSu},
 if the radius $\rho$ is sufficiently large, then 
\bel{redd}\S(\mu_k^\flat) - c\I^\alpha(\mu_k^\flat)~\geq~\S(\mu_k) - c\I^\alpha(\mu_k)\eeq
 for all $k\geq 1$. Here $\S$ and $\I^\alpha$ are the functionals defined at 
 (\ref{S8})-(\ref{TC}). According to (\ref{redd}), we can replace the measure $\mu_k$ with $\mu_k^\flat$
 without decreasing  the objective functional.
 
 Without loss of generality we can thus choose $\ell>0$ sufficiently large and assume that 
	$$u_k(s) = 0 \qquad \forall s > \ell, \quad k \ge 1.$$
In turn, since $\S(\mu_k)-c\I^\alpha(\mu_k)\geq 0$, we obtain the uniform bound
\bel{iab}
\I^\alpha(\mu_k)~\leq~\kappa_1~\doteq~{1\over c} \S(\mu_k)~\leq~{\ell\over c}\,.\eeq
\v
{\bf 2.} In this step we show that the measures $\mu_k$ can be taken with uniformly bounded mass.
Consider a measure $\mu_k$ for which (\ref{iab}) holds.
By  (\ref{TC}), for every $r\in [0,\ell]$  one has 
$$\I^\alpha(\mu_k)~\geq~r\cdot \left(\int_r^\ell u_k(t)\, dt\right)^\alpha.$$
In view of (\ref{iab}), this implies
\bel{rki} \int_{r}^{\ell}u_k(s)\,ds~ \le ~\left( {\kappa_1 \over r} \right)^{1/\alpha}.\eeq
It thus remains to prove that, in our maximizing sequence, the functions $u_k$ can be replaced with 
functions $\tilde u_k$ having a uniformly bounded integral over  $[0,r]$, 
for some fixed $r>0$.

Toward this goal we fix $0 < \ve < \beta < 1$, and, for $j\geq 1$, we  define  $r_j = 2^{-j}$,  and the 
interval
$V_{j} =\, ]r_{j+1},r_j]$.   Given $u=u_k$, 
if  $\int_{V_j}u(s)\,ds > r_j^\ve$,
we introduce the functions
\bel{ujj} u_j(s) \,\doteq\, \chi_{\strut V_j}(s) u(s),\qquad\qquad \tilde{u}_j (s)\,\doteq
\, \min\{ u_j(s), c_j\},%2 r_j^{\beta-1} \chi_{\strut V_j}(s),
\eeq
choosing the constant $c_j\geq 2 r_j^{\beta-1}$ so that	
\bel{irj} \int_{V_j} \tilde{u}_j(s)\,ds ~=~ r_j^\beta.\eeq
We then let $\mu_j= u_u\mu$ and $\tilde{\mu}_j= \tilde u_j\mu$ be the measures 
supported on $V_j$, corresponding to these densities.

For a fixed integer $j^*$, whose precise value will be chosen later,
consider the set of indices
\bel{Jdef} J~ \doteq ~\left\{ j \ge j^* \:\bigg|\: \int_{V_j}u(s)\,ds > r_j^\ve  \right\} \eeq
	and the modified density
	\bel{mum} \tilde{u}(s) \doteq u(s) + \sum_{j \in J} ( \tilde{u}_j(s) - u_j(s) ). \eeq
Moreover, call $\tilde \mu$ the  measure obtained by replacing $u$ with $\tilde u$ in (\ref{m2}).
By (\ref{rki}) and (\ref{ujj}) the total mass of $\tilde \mu$ is bounded. Indeed
\bel{utu}\tilde\mu(\R^2)~=~
 \int_{r_{j^*}}^{\ell}\tilde{u}(s)\,ds + \int_{0}^{r_{j^*}}\tilde{u}(s)\,ds ~\le~
  \left( {\kappa_1 \over r_{j^*}} \right)^{1/\alpha} +  \sum_{j \ge j^*} r_j^\ve ~\le~ \left( {\kappa_1 \over r_{j^*}} \right)^{1/\alpha} + \sum_{j \geq 1 } 2^{-j\ve} ~%\doteq M 
  <~ +\infty. \eeq

 We now claim that
 \bel{mtm}\S(\tilde \mu) - c\I^\alpha(\tilde \mu)~\geq~\S(\mu) - c\I^\alpha(\mu).\eeq
Toward a proof of (\ref{mtm}), we %fix any $j\in J$ and 
estimate
\bel{Sbd6}\bega{rl} \S(\mu) - \S(\tilde \mu)&\leq \ds~\sum_{j\in J} \bigg(
	 \int_{V_j} I(y(t)) \cos(\theta(t)-\theta_0)\,dt   \\ [4mm] \ds
	&\ds \qquad\qquad - \int_{V_j} I(y(t)) \left(1 - \exp\Big\{ -{\tilde{u}_j(t) \over \cos(\theta(t)-\theta_0)} \Big\} \right) \cos(\theta(t)-\theta_0)\,dt\bigg) \\ [4mm]
	&\ds\leq~ \sum_{j\in J}\int_{r_{j+1}}^{r_{j}} \exp\bigl\{ -\tilde{u}_j(t) \bigr\} dt ~
	\leq~ \sum_{j\in J}r_{j+1} \exp\left\{-2 r_j^{\beta-1} \right\}. \enda\eeq
To estimate the difference in the irrigation cost, we first observe that the inequality
$$\left( \int_r^\ell u(t)\, dt\right)^\alpha~\leq~{1\over r} \I^\alpha(\mu)~=~{\kappa_1\over r}$$
implies
\bel{uib}
\left( \int_r^\ell u(t)\, dt\right)^{\alpha-1}~\geq~\left({\kappa_1\over r}\right)^{\alpha-1\over\alpha}\,.\eeq
Since $\tilde u(s)\leq u(s)$ for every $s\in [0,\ell]$, using (\ref{uib}) we now obtain
\bel{IIT}\bega{rl} \I^\alpha(\mu)-\I^\alpha(\tilde \mu) &\ds
=~\int_0^1{d\over d\lambda} \I^\alpha\bigl(\lambda\mu + (1-\lambda) 
\tilde \mu\bigr)\, d\lambda\\[4mm]
&=~\ds \int_0^1\int_0^\ell {d\over d\lambda}\left( \int_s^\ell [\lambda u(t) + (1-\lambda) 
\tilde u(t)]\, dt\right)^\alpha ds\, d\lambda\\[4mm]
&=~\ds\int_0^1 \int_0^\ell\left\{ \alpha \left( \int_s^\ell [\lambda u(t) + (1-\lambda) 
\tilde u(t)]\, dt\right)^{\alpha-1} \int_s^\ell [u(t)-\tilde u(t)]\, dt \right\}ds\, d\lambda\\[4mm]
&\geq~\ds \int_0^\ell\left\{ \alpha \left( \int_s^\ell u(t)\, dt\right)^{\alpha-1} \int_s^\ell [u(t)-\tilde u(t)]\, dt \right\}ds\\[4mm]
&\ds\geq~\sum_{j\in J}
\int_{r_{j+2}}^{r_{j+1}} \left[ \alpha \left(  \int_s^{\ell} u(t)\,dt\right)^{\alpha-1} \int_{r_{j+1}}^{r_j} (u_j(t) - \tilde{u}_j(t) ) 
\,dt \right] ds \\[4mm]
	&\ds\geq~ \sum_{j\in J} \alpha\left( {\kappa_1\over r_{j+2}}\right)^{\alpha-1\over \alpha}\cdot (r_j^\ve - r_j^\beta)
\cdot r_{j+2}\\[4mm]	
	&\ds=~\sum_{j\in J}  \kappa_2  r_j^{1/\alpha} (r_j^\ve - r_j^\beta),\enda
\eeq
where $ \kappa_2=\alpha(4\kappa_1)^{\alpha-1\over\alpha}$.   Combining (\ref{Sbd6}) with 
(\ref{IIT}) we obtain
\bel{SS2} c[\I^\alpha( \mu) - \I^\alpha(\tilde\mu)]- [\S( \mu) -\S(\tilde\mu)]    ~\geq~
\sum_{j\in J} \Big(c\kappa_2  r_j^{1/\alpha} (r_j^\ve - r_j^\beta)- r_{j+1} \exp\left\{-2 r_j^{\beta-1} \right\}\Big).
\eeq
By choosing the integer $j^*$ large enough in (\ref{Jdef}), for $j\geq j^*$ all terms in the summation 
on the right hand side of (\ref{SS2}) are $\geq 0$.  This implies (\ref{mtm}).
\v
{\bf 3.} By the two previous steps, w.l.o.g.~we can assume that the measures $\mu_k$ have uniformly bounded support and  uniformly bounded total mass.  Otherwise,  
we can replace the sequence $(u_k)_{k\geq 1}$ with a new maximizing sequence
$(\tilde u_k)_{k\geq 1}$  having these properties.

By taking a subsequence, we can thus assume the weak convergence $\mu_k\wto\ov \mu$.
The upper semicontinuity of the functional $\S$, proved in \cite{BS}, yields
\bel{SUC}
\S(\ov \mu)~\geq~\limsup_{k\to\infty} ~\S(\mu_k).\eeq
In addition, since all maps $s\mapsto \gamma_k(s)$ are 1-Lipschitz, by taking a further 
subsequence we can assume the convergence
\bel{gkg}\gamma_k(s)~\to~\ov\gamma(s)\eeq
for some limit function $\ov\gamma$, uniformly for $s\in [0,\ell]$.

Since each measure $\mu_k$ is supported on $\gamma_k$,
the weak limit $\ov \mu$ is a measure supported on the 
curve $\ov\gamma$.   %$\bar \gamma=\{(x,y)\,;~~x = \bar x(y)\,;~~y\in [0,\ell]\}$.
\v
{\bf 4.}  Since $\theta_k(s)\in [\theta_0,\pi/2]$,
we can re-parameterize each stem $\gamma_k$ in terms of the vertical variable
$$y_k(s)~=~\int_0^s \sin\theta_k(s)\, ds.$$
Calling $s=s_k(y)$ the inverse function, we thus obtain a maximizing sequence
of couples
$$y~\mapsto~(\Hat \theta_k(y), \Hat u_k(y))~\doteq~
\Big(\theta_k(s_k(y)),\, u_k(s_k(y))\Big),\qquad\qquad y\in [0, h_k]\,.$$
Moreover, the stem $\gamma_k$ can be described as the graph of the Lipschitz function
$$x~=~x_k(y)~=~\int_0^{s_k(y)} \cos\theta_k(s)\, ds.$$
Since all functions $x_k(\cdot)$ satisfy $x_k(0)=0$ and are non-decreasing,
uniformly continuous with Lipschitz constant
$L=\cos\theta_0/\sin \theta_0$,
by possibly extracting a further subsequence, we obtain the convergence $h_k\to \bar h$ and 
%the uniform convergence
$x_k(\cdot) \to \bar x(\cdot)$.  Here $\bar x:[0, \bar h]\mapsto \R$ is a nondecreasing  continuous function
with Lipschitz constant $L$, such that $\bar x(0)=0$.     More precisely, the convergence $x_k\to \bar x$
is uniform on every compact subinterval $[0,h]$ with $h<\bar h$.
\v
{\bf 5.} We claim that  the irrigation cost of $\ov\mu$ is no greater that the
lim-inf of the irrigation costs for $\mu_k$.   Let $\sigma\mapsto \gamma(\sigma)$ 
be an arc-length
parameterization of $\ov \gamma$. Since $s\mapsto \ov \gamma(s)$ if 1-Lipschitz, 
one has $d\sigma/ds\leq 1$.
We now compute
\bel{ILC}\bega{rl}\I^\alpha(\ov\mu)&=~\ds
\int_0^{\sigma(\ell)} \left(  \int_\sigma^{\sigma(\ell)} \ov u(t)\,dt\right)^\alpha
d\sigma~=~\int_0^{\sigma(\ell)} \left(  \lim_{k\to\infty} \int_s^\ell  u_k(t)\,dt\right)^\alpha
d\sigma(s)\\[4mm]
&\ds\leq~\lim_{k\to\infty} \int_0^\ell \left(\int_s^\ell  u_k(t)\,dt\right)^\alpha
ds~=~\lim_{k\to\infty} \I^\alpha(\mu_k).\enda\eeq
\v
{\bf 6.} Combining (\ref{SUC}) with (\ref{ILC}) we conclude that the measure $\ov \mu$, supported on the stem $\ov\gamma$, is optimal.

Let $\bar u$ be the density of the absolutely continuous part of $\ov\mu$  w.r.t.~the arc-length measure on $\bar \gamma$, and call $\mu^*$ the measure that has density $\bar u$ w.r.t.~arc-length measure.
Since $\S(\mu^*)=\S(\ov\mu)$, it follows that $\mu^*=\ov\mu$. Otherwise $\I^\alpha(\mu^*)<\I^\alpha(\ov\mu)$
and $\ov\mu$ is not optimal.    
This argument shows that the optimal measure $\ov \mu$ is absolutely continuous w.r.t.~the arc-length 
measure on $\ov\gamma$.

Calling $\sigma\mapsto \gamma(\sigma)$ the arc-length parameterization of $\ov\gamma$,
the optimal solution to {\bf (OP2)}  is now provided by 
$\sigma\mapsto(\ov \theta(\sigma), \bar u(\sigma))$, where
$\ov\theta$ is the orientation of the tangent vector:
$${d\over d\sigma} \ov\gamma(\sigma)~=~\bigl(\cos \ov \theta(\sigma), \, \sin\ov \theta(\sigma)
\bigr).$$
\endproof

\subsection{Necessary conditions for optimality.}

Let $t\mapsto (\theta^*(t), u^*(t))$ be an optimal solution to
 the problem {\bf (OP2)}.   % stated at (\ref{max33})-(\ref{max44}).
The necessary conditions for optimality \cite{BP, Cesari, FR} yield the existence of dual variables $p,q$ satisfying
\bel{pq1}
\left\{\bega{rl} \dot p&=~- I'(y)\, G(\theta,u),\\[3mm]
\dot q&=~ c \alpha \,z^{\alpha-1},
\enda\right.\qquad\qquad 
\left\{ \bega{cl} p(+\infty)&=~0,\\[3mm]
q(0)&=~0,\enda\right.\eeq
and such that the maximality condition
\bel{MC}
(\theta^*(t), u^*(t))~=~\hbox{arg}\!\!\!\!\!\!\!\max_{\theta\in [0,\pi], ~u\geq 0}~
\Big\{ p(t)\, \sin\theta - q(t) u + I(y(t))\, G(\theta,u)- c z^\alpha
\Big\}.\eeq
We recall that 
$G(\theta,u)$ 
%~\doteq~\left( 1 - \exp\left\{-{u\over \cos(\theta-\theta_0)}
%\right\}\right)\cos(\theta-\theta_0)$$
is the function defined at (\ref{gtu}).
An intuitive interpretation of the quantities on the right had side of (\ref{MC}) goes as follows:
\begi
\item $p(t)$ is the rate of increase in the gathered sunlight, if the upper portion of stem $\{\gamma(s)\,;~~s>t\}$ is
raised higher.
\item $q(t)$ is the rate at which the irrigation  cost increases, adding mass at the point $\gamma(t)$.
\item  $I(y(t))\, G(\theta,u)$ is the sunlight captured by the leaves at the point $\gamma(t)$.
\endi

\section{Uniqueness of the optimal stem configuration}  \label{s:55}
\setcounter{equation}{0} 
Aim of this section is to  show that,
if the light intensity $I(y)$ remains sufficiently close to 1 for all $y\geq 0$, then the shape of the optimal stem is uniquely determined.    This models a case where the density of external vegetation is small.

\begin{theorem}\label{t:61}
Let $h\mapsto I(h)\in [0,1]$ be a non-decreasing, absolutely continuous
 function which satisfies
\bel{Iprime} I'(y)\,\leq \,C y^{-\beta} \qquad\hbox{for a.e.}~ y>0,\eeq
for some constants $C>0$ and $0<\beta <1$.
If
\bel{sveg} I(0)\,\geq \,1-\delta\eeq
for some $\delta>0$ sufficiently small,
then the optimal solution
to {\bf (OP2)} is unique.
\end{theorem}

{\bf Proof.}
We will show that the necessary conditions for optimality have a unique solution.
This will be achieved in several steps.
\v
{\bf 1.} 
Given $I,p,q$, define the functions $\Theta, U$ by setting
\bel{max41} \Big(\Theta(I,p,q),\, U(I,p,q)\Big)~\doteq~
\argmax_{\theta\in [0,\pi], ~u\geq 0}
\Big\{ p\cdot \sin\theta - q\, u + I\cdot G(\theta,u)- c z^\alpha\Big\}.\eeq
We recall that $G$ is the function defined at (\ref{gtu}).
Notice that one can write
$$
G(\theta,u) \,=\, u \,\Tilde G\left(\frac{\cos{\left(\theta-\theta_0\right)}}{u}\right)$$
with
\bel{TGprop}\Tilde G(x) \,\doteq \,\left(1-\exp\left\{-\frac{1}{x}\right\}\right)x~>~0,
\quad \qquad  \Tilde G' (x)\,\le\, 1,\qquad \Tilde G''(x)\, \le\, 0,\qquad\forall x>0.\eeq

Denote by
\bel{PSI}{\cal H}(\theta,u)~\doteq~ p\cdot \sin\theta - q\, u + I(y)\, G(\theta,u)- c z^\alpha\eeq
 the quantity  to be maximized in (\ref{max41}).   Differentiating ${\cal H}$
 w.r.t.~$\theta$ and imposing that
the derivative is zero, we obtain
\bel{Theta1}\bega{rl}\ds
\frac{p}{I} &\ds=\, -\frac{ G_\theta(\theta,u)}{\cos{\theta}} \\[4mm]
&= ~\ds\frac{\sin{(\theta-\theta_0)}}{\cos{\theta}} \left[1 - \exp\left\{ -\frac{u}{\cos{(\theta-\theta_0)}} \right\} - \frac{u}{\cos{(\theta-\theta_0)}} \exp\left\{ -\frac{u}{\cos{(\theta-\theta_0)}} \right\} \right].\enda\eeq
Similarly, differentiating w.r.t.~$u$  we find %also find when the partial derivative with respect to $u$ is zero, i.e.\
\[
-q + I G_u(\theta,u) ~= ~-q + I \exp\left\{ -\frac{u}{\cos{(\theta-\theta_0)}} \right\}~ = ~0.
\]
This yields
\bel{u7}
u~ = ~-\ln\left(\frac{q}{I}\right)\,\cos{(\theta-\theta_0)}.
\eeq
A lengthy but elementary computation shows that the Hessian matrix
of second derivatives of ${\cal H}$ w.r.t.~$\theta,u$ is negative definite, and the critical point 
is indeed the  point where the global maximum is attained. By (\ref{u7}) it follows
\begin{equation}\label{U}
U(I,p,q) ~= ~-\ln\left(\frac{q}{I}\right)\,\cos\bigl(\Theta(I,p,q) - \theta_0 \bigr).
\end{equation}
Inserting (\ref{U}) in  (\ref{Theta1}) and using the identity
\[
\frac{\sin{(\theta-\theta_0)}}{\cos{\theta}}~ =~ \cos{\theta_0} \tan{\theta} - \sin{\theta_0}
\]
we obtain
\begin{equation}\label{Theta2}
\Theta(I,p,q) ~=~ \arctan\left( \tan{\theta_0} + \frac{ {1\over \cos\theta_0}\,\frac{p}{I}}{1-\frac{q}{I} +\frac{q}{I} \ln{\left(\frac{q}{I}\right)} } \right)
%~=~\arctan\left( \tan \theta_0 + \frac{ p/I}{ \cos\theta_0\cdot w(q/I)} \right),
\end{equation}
Introducing the function
	\begin{equation}
	w(I,p,q)~ \doteq ~\frac{p/I}{1-\frac{q}{I} +\frac{q}{I} \ln{\left(\frac{q}{I}\right)} } ,
	\label{w}
	\end{equation}
by  (\ref{Theta2}) one has the identities
\bel{trig1}\left\{\bega{rl}
	\sin\bigl(\Theta(I,p,q)\bigr)&=~\ds \frac{\sin{\theta_0} + w}{ \sqrt{ \cos^2{\theta_0} + (w + \sin{\theta_0})^2 } }, \\[4mm]
	\cos\bigl(\Theta(I,p,q)-\theta_0\bigr) &=\ds \frac{1 + w\sin{\theta_0} }{ \sqrt{ \cos^2{\theta_0} + (w + \sin{\theta_0})^2 } },
	\enda\right.\eeq
	Note that $w \ge 0$, because $p, q, I \ge 0$. In turn, 
from (\ref{trig1}) it follows
\bel{trig2}\left\{\bega{rl}
	\cos\bigl(\Theta(I,p,q)\bigr)&=~\ds\frac{\cos{\theta_0}}{ \sqrt{ \cos^2{\theta_0} + (w + \sin{\theta_0})^2 } }, \\[4mm]
	\sin\left(\Theta(I,p,q)-\theta_0\right) &=\ds~ \frac{w\cos{\theta_0} }{ \sqrt{ \cos^2{\theta_0} + (w + \sin{\theta_0})^2 } }.
	\enda\right.\eeq

%where 
%\bel{wdef} w(\zeta)~\doteq~1-\zeta+\zeta\,\ln |\zeta|.\eeq
\v
{\bf 2.}
The necessary conditions for the optimality of a solution to {\bf (OP2)}
yield the boundary value problem
\bel{BVP}\left\{\bega{rl}
\dot{y}(t) &=~ \sin\Theta, \\[3mm]
\dot{z}(t) &= ~-U,\\[3mm]
\dot{p}(t) &= ~-I'(y) G\bigl( \Theta, U \bigr), \\[3mm]
\dot{q}(t) &= ~c \alpha z^{\alpha-1},
\enda \right.\qquad\qquad 
\left\{
\bega{rl}
y(0) &=~ 0,\\[3mm]
z(T) &=~ 0,\\[3mm]
p(T) &= ~0,\\[3mm]
q(T) &= ~I(y(T)),\\[3mm]
q(0) &=~ 0.
\enda \right.\eeq
Here  $[0,T[\,$ is the interval where $u>0$, while
\bel{TU}
\Theta\,=\,\Theta(I(y), p,q),\qquad U\,=\, U(I(y), p,q)\eeq
are the functions introduced at (\ref{max41}), or more explicitly at 
(\ref{U})-(\ref{Theta2}).  Notice that the length $T$ of the stem is a quantity to be determined, 
using the boundary conditions in (\ref{BVP}).
\v
{\bf 3.} 
Since the control system (\ref{max44}) and the running cost (\ref{max33}) 
do not depend explicitly on time, the Hamiltonian function
\bel{ham}H(y,z,p,q) ~ \doteq~\max_{\theta\in [0,\pi], ~u\geq 0}
\Big\{ p\cdot \sin\theta - q\, u + I(y)\, G(\theta,u)- c z^\alpha\Big\}\eeq 
is constant along trajectories of (\ref{BVP}).
Observing that the terminal conditions in (\ref{BVP})  imply 
$H(y(T),z(T),p(T),q(T)) = 0$, one has the first integral
\bel{H0}H(y(t),z(t),p(t),q(t))~=~0\qquad\qquad\forall t\in [0,T].\eeq
This yields
$$\bega{rl}
	0 &= \ds~p \sin{\Theta} + \left[I(y) - q + q \ln\left(\frac{q}{I(y)}\right)\right] \cos\left(\Theta-\theta_0\right) - c z^\alpha \\[4mm]
	&=~\ds \frac{p \left[\sin{\theta_0} + w\right] + \left[I(y) - q + q \ln\left(\frac{q}{I(y)}\right)\right] \left[1 + w \sin{\theta_0} \right]}{\sqrt{ \cos^2{\theta_0} + (w + \sin{\theta_0})^2 }} - c z^\alpha\\[4mm]
	&=~\ds I(y) \left[1 - \frac{q}{I(y)} + \frac{q}{I(y)} \ln\left(\frac{q}{I(y)}\right)\right] \sqrt{ \cos^2{\theta_0} + (w + \sin{\theta_0})^2 } - c z^\alpha.\enda
	$$
We can use this identity to express $z$ as a function of the other variables:
\bel{z}\bega{rl}
	z\bigl(I(y),p,q\bigr) &= ~\ds\bigg\{ \frac{I(y)}{c} \left[1 - \frac{q}{I(y)} + \frac{q}{I(y)} \ln\left(\frac{q}{I(y)}\right)\right] \sqrt{ \cos^2{\theta_0} + (w + \sin{\theta_0})^2 } \bigg\}^{1/\alpha} \\[4mm]
	&=~\ds c^{-1/\alpha} \bigg\{ \left( \left[I(y) - q + q \ln\Big(\frac{q}{I(y)}\Big)\right] \cos{\theta_0}\right)^2   \\[4mm]
	&\qquad\qquad \qquad
	+ \left(p + \left[I(y) - q + q \ln\left(\frac{q}{I(y)}\right)\right] \sin{\theta_0} \right)^2 \bigg\}^{1/2\alpha}.
	\enda
\eeq

\begin{figure}[ht]
\centerline{\hbox{\includegraphics[width=15cm]
{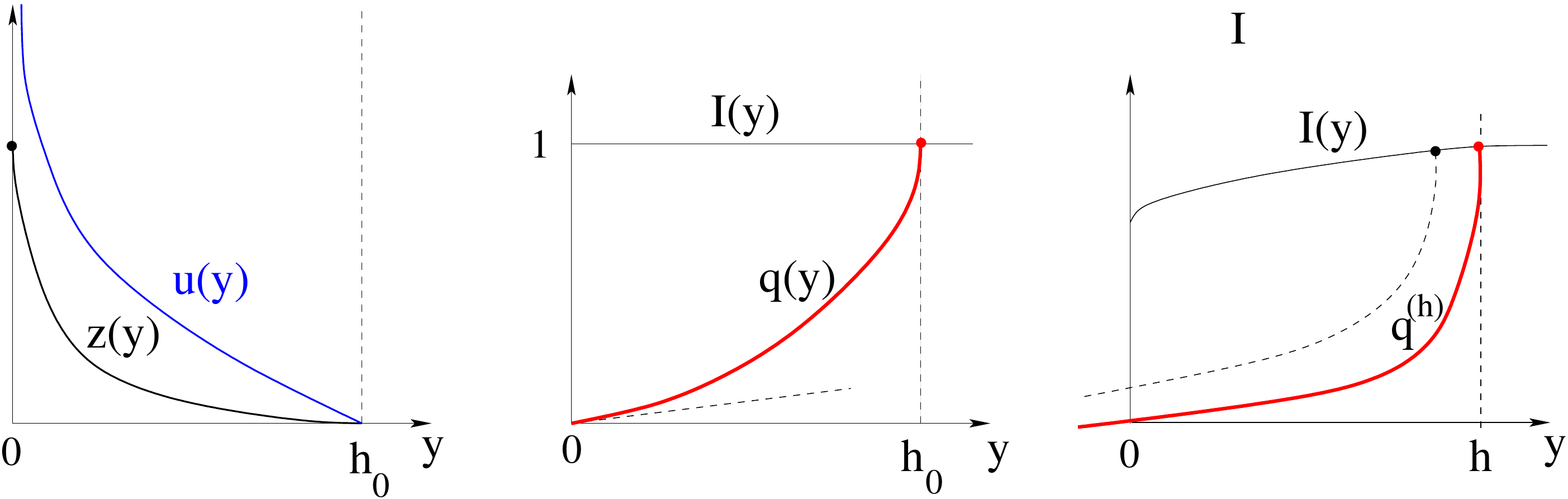}}}
\caption{\small  Left and center: sketch of the solution of  the system (\ref{pq1}) 
in the case where $I(y)\equiv 1$.  Left: the graphs of the functions $z$ in (\ref{z10}) and $u= -\ln z$.
Center: the graph of the function $q$ at (\ref{qy}).   The figure on the right shows the case where
$I(\cdot)$ is not constant.  As before,  $h$ must be determined so that $q(0,h)=0$. }
\label{f:pc1}
\end{figure}

\v{\bf 4.} 
Since $I$ is given as a function of the height $y$, it is convenient to rewrite the equations 
(\ref{BVP}) using $y$ as an independent variable.  
Using  the identity (\ref{z}), we obtain a system
 of two equations for the variables $p,q$:
\bel{dpdy}	\bega{rl} \ds{d\over dy} p(y)
 &=\ds~ -I'(y) \left[1-\frac{q(y)}{I(y)} \right] \frac{\cos\left(\Theta\bigl(I(y),p(y),q(y)\bigr)
  - \theta_0\right)}
 {\sin \Theta\bigl(I(y),p(y),q(y)\bigr)} \\[4mm]
	&=~\ds -I'(y) \left[1-\frac{q(y)}{I(y)} \right] \frac{1 + w \sin{\theta_0}}
	{w + \sin{\theta_0}} \\[4mm]
	&\doteq~\ds -I'(y) \,f_1\bigl(I(y),p(y),q(y)\bigr) ,
	\enda\eeq
\bel{dqdy}	\bega{rl} \ds{d\over dy} q(y)
 &=~\ds \frac{c \alpha \bigl[z\bigl(I(y),p(y),q(y)\bigr)\bigr]^{\alpha-1}}
 {\sin\Theta\bigl(I(y),p(y),q(y)\bigr)} \\[4mm]
	&=~\ds \frac{\alpha c^{1/\alpha}}{w + \sin{\theta_0}} \left[ \cos^2{\theta_0} + (\sin{\theta_0} + w)^2 \right]^{1-\frac{1}{2\alpha}} \\[4mm]
	&\ds\qquad\qquad\times \left[ I(y) \left( 1-\frac{q}{I(y)} + \frac{q}{I(y)}\ln\left(\frac{q}{I(y)}\right) \right) \right]^{1-\frac{1}{\alpha}} \\[4mm]
	&\doteq~ f_2\bigl(I(y),p(y),q(y)\bigr),
\enda\eeq
where $w = w(I,p,q)$ is the function introduced at (\ref{w}).
Note that under our assumptions, $f_1$ remains  bounded, while $f_2$ diverges as $q(y) \to I(y)$.
The system (\ref{BVP}) can now be equivalently formulated as
 \bel{bvp2}
\left\{\bega{rl}\ds p'(y)&=\,-I'(y)\, f_1\bigl(I(y),p,q\bigr), \\[3mm]
 q'(y) & =~f_2\bigl(I(y),p,q\bigr),
\enda \right.\qquad\qquad \left\{
\bega{rl}
p(h) &= ~0,\\[3mm]
q(h) &= ~I(h),
\enda \right.\qquad\quad
q(0) \,=\,0.\eeq
\v
{\bf 5.}
To prove uniqueness of the solution to the boundary value problem (\ref{BVP}), 
it thus suffices to prove the following (see Fig.~\ref{f:pc1}, right).
\begi
\item[{\bf (U)}]{\it
Call 
\bel{pqyh}y~\mapsto~\bigl( p(y,h), \, q(y,h)\bigr) \eeq
the solution to the system} (\ref{bvp2}), {\it with the two terminal conditions
given at $y=h$.    Then there is a unique choice of $h>0$ which satisfies also
the third boundary condition}
\bel{q0}q(0,h)~=~0.\eeq
\endi
To make the argument more clear, the uniqueness property {\bf (U)} will be proved in two steps.
\begi
\item[(i)] 
When $I(y)\equiv 1$,
the map 
\bel{hq}
h~\mapsto q(0,h)\eeq
is strictly decreasing, hence it vanishes at a unique point $h_0$.
\item[(ii)] For all functions $I(\cdot)$ sufficiently close to the constant map $\equiv 1$, the map
(\ref{hq}) is strictly decreasing in a neighborhood of $h_0$.
\endi
\v
In the case $I(y)\equiv 1$, recalling (\ref{Theta2}) we obtain (see Fig.~\ref{f:pc1})
$$I'(y)~=~0,\qquad p(y,h)~=~0,\qquad \Theta(I, 0,q)~=~\theta_0,\qquad G(\theta_0, U)~=~1-e^{-U},$$
$$U(1,0,q)~=~\argmax_u\,\bigl\{ -qu + G(\theta_0,U)\bigr\}~=~\argmax_u\{ -qu +1-e^{-u}\}~=~-\ln q,$$
The system (\ref{BVP}) can now be written as
\bel{ODE20}\left\{\bega{rl}  p'(y)& = ~0, \\[3mm]
 q'(y)&=~\ds{c \alpha z^{\alpha-1}\over\sin\theta_0}\,,\\[3mm]
 z'(y)&=~ \ds {\ln q\over \sin\theta_0}.
\enda \right.
\qquad \qquad
\left\{
\bega{rl}
p(h) &= ~0,\\[3mm]
q(h) &= ~1,\\[3mm]
z(h)&=~0,
\enda \right.\qquad q(0) ~=~ 0.\eeq
{}From (\ref{ODE20}) it follows   $p(y)~\equiv~0$, while
$$ {dz\over dq}~=~{\ln q\over c\alpha z^{\alpha-1}}.$$
Integrating the above ODE with terminal conditions $q=1$, $z=0$, one obtains
\bel{z10}
z~=~c^{-1/\alpha} \Big[1+ q\ln q-q\Big]^{1/\alpha}.\eeq
The second equation in (\ref{ODE20}) thus becomes
\bel{qy}
q'(y)~=~{\alpha c^{1/\alpha}\over\sin\theta_0}   \Big[1+ q\ln |q|-q\Big]^{\alpha-1\over\alpha}.\eeq
Notice that here  the right hand side is strictly positive for all $q\in \,]-1,1[\,$.   Of course,  only positive values of 
$q$ are relevant for the optimization problem, but for the analysis it is convenient
to extend the definition also to negative values of $q$.
The solution of (\ref{qy}) with terminal condition $q(h)=1$ is implicitly determined by
\bel{qys} h-y~=~
{\sin\theta_0\over \alpha c^{1/\alpha}}\int_{q(y)}^1\Big[ 1+s\ln |s|-s\Big]^{1-\alpha\over\alpha}\, ds
\,.\eeq
The map $h\mapsto q(0,h)$  thus vanishes at the unique point
\bel{h0}h_0~=~
{\sin\theta_0\over \alpha c^{1/\alpha}}\int_0^1\Big[ 1+s\ln |s|-s\Big]^{1-\alpha\over\alpha}\, ds.\eeq
As expected, the height $h_0$ of the optimal stem decreases as we increase the constant $c$, 
in the transportation cost.
A straightforward computation yields
\bel{dq0} {\partial \over \partial h} q(0,h)~=~-
{ \alpha c^{1/\alpha}\over \sin\theta_0} \Big[ 1+q(0,h)\ln |q(0,h)|-q(0,h)\Big]^{1-\alpha\over\alpha}.
\eeq
In particular, at $h=h_0$ we have  $q^{(h_0)}(0)=0$ and hence
\bel{dq00} {d\over dh} q(0,h)\bigg|_{h=h_0}~=~-
{ \alpha c^{1/\alpha}\over \sin\theta_0}~<~0.
\eeq

\v
{\bf 6.} We will show that a strict inequality as in (\ref{dq00}) remains valid for a more general 
function $I(\cdot)$, provided that the  assumptions (\ref{Iprime})-(\ref{sveg}) hold.

Toward this goal, we need to determine how $p$ and $q$ vary w.r.t.~the parameter $h$.
Denoting by 
\bel{PQdef}P(y) ~\doteq~ \pdiff{p(y,h)}{h}, \qquad Q(y)~ \doteq~ \pdiff{q(y,h)}{h}\eeq
their partial derivatives, by (\ref{bvp2}) one obtains the linear system
	\bel{PQ'}
	\begin{pmatrix}
	P(y) \\ Q(y)
	\end{pmatrix}' ~= ~\begin{pmatrix}
	-I'(y) f_{1,p} &&- I'(y) f_{1,q} \\ f_{2,p} && f_{2,q}
	\end{pmatrix} \begin{pmatrix}
	P(y) \\ Q(y)
	\end{pmatrix}.
	\eeq
The boundary conditions at $y=h$ require some careful consideration.
As $y\to h-$, we expect $f_2(I(y), p(y), q(y))\to +\infty$ and $Q(y)\to -\infty$.
To cope with this singularity we  introduce the new variable
\bel{TQdef}\Tilde Q(y)~ \doteq~ {Q(y) \over f_2\bigl(I(y),p(y),q(y)\bigr)}\,.\eeq
	The system (\ref{PQ'}), together with the new boundary conditions for $P,\Tilde Q$
	can now be written as
\bel{dh1}\left\{\bega{rl}
	P'(y) &=\,\ds- I'(y) \left[ f_{1,p} P + f_{1,q} f_2 \Tilde{Q} \right], \\[3mm]
	\Tilde{Q}'(y) &=~\ds \frac{f_{2,p}}{f_2} P - \frac{I'(y) [f_{2,I} - f_{2,p} f_1]}{f_2} \Tilde{Q},
	\enda\right.\qquad\qquad \left\{ \bega{rl} P(h)&=~0,\\[3mm]
	\Tilde Q(h)&=~-1.\enda\right.\eeq
	To analyze this system we must compute the partial derivatives of $f_1$ and $f_2$.
 From the definition (\ref{w}) it follows
 \bel{partw}\pdiff{w}{I} ~=~ \frac{w^2}{p} \left[1 - \frac{q}{I}\right],\qquad
	\pdiff{w}{p} ~=~ \frac{w}{p}, \qquad
	\pdiff{w}{q} ~=\, -\frac{w^2}{p} \ln\left(\frac{q}{I}\right).\eeq
Using  (\ref{partw}), from \eqref{dpdy}, \eqref{dqdy} we obtain
	\bel{fpart}\left\{
	\bega{rl}
	f_{1,p}\bigl(I(y),p,q\bigr) &=\ds ~\frac{ 1 - \frac{q}{I(y)} }{I(y) \tan^2{\Theta} \left[ 1 - \frac{q}{I(y)} + \frac{q}{I(y)}\ln\left(\frac{q}{I(y)}\right) \right] }, \\[4mm]
	f_{1,q}\bigl(I(y),p,q\bigr)  &=\ds ~\frac{1}{I(y)} \frac{\cos\left(\Theta-\theta_0\right)}{\sin\Theta} - \frac{ \sin\left(\Theta-\theta_0\right) \cos\Theta \left[1 - \frac{q}{I(y)}\right] \ln\left(\frac{q}{I}\right) }{I(y) \sin^2{\Theta} \left[ 1 - \frac{q}{I(y)} + \frac{q}{I(y)}\ln\left(\frac{q}{I(y)}\right) \right] },
\\[4mm]
	f_{2,p}\bigl(I(y),p,q\bigr)  &\ds= ~
	-\left[1 + \frac{\alpha}{\sin^2\Theta} - 2\alpha \right] \frac{1}{z\bigl(I(y),p,q\bigr) },\\[4mm]
	f_{2,q}\bigl(I(y),p,q\bigr) ) &= \ds ~-\left[ \frac{(1-\alpha) \sin{\theta_0}}{\sin^2\Theta} - \frac{\sin\left(\Theta-\theta_0\right)}{\cos\Theta}\left(1 + \frac{\alpha}{\sin^2\Theta} -2\alpha \right) \right] \frac{ \ln\left(\frac{q}{I(y)}\right) }{z\bigl(I(y),p,q\bigr) } ,\\[4mm]
	f_{2,I}\bigl(I(y),p,q\bigr)  &= \ds ~ -\left[\frac{(1-\alpha) \sin{\theta_0}}{\sin^2\Theta} + \frac{\sin\left(\Theta-\theta_0\right)}{\cos\Theta}\left(1 + \frac{\alpha}{\sin^2\Theta} -2\alpha \right)  \right] \frac{1-\frac{q}{I(y)}}{z\bigl(I(y),p,q\bigr) } .
	\enda\right.\eeq
At this stage, the strategy of the proof is straightforward.
When $I'(y)\equiv 0$, the solution to (\ref{dh1}) is trivially given by $P(y)\equiv 0$,
$\Tilde Q(y)\equiv -1$.    This implies
$${\partial\over\partial h}q(0,h)~=~\Tilde Q(0)\cdot f_2(I(0), p(0), q(0))~<~0.$$
We need to show that the same strict inequality holds when $\delta>0$ in (\ref{sveg}) 
is small enough. 
Notice that, if the right hand sides of the equations
in (\ref{dh1}) were bounded, 
letting  $\|I'\|_{\L^\infty}\to 0$ a continuity argument would imply the uniform convergence
$P(y)\to 0$  and $\Tilde Q(y)\to -1$.  
The same conclusion can be achieved provided that 
the right hand sides in (\ref{dh1}) are uniformly integrable. 
This is precisely what will be proved in the next two steps, relying on the 
identities (\ref{fpart}).
\v
{\bf 7.} In this step we prove an inequality of the form
\bel{TUB}
0~<~\theta_0~\leq~\Theta(I,p,q)~\leq ~\theta^+~<~{\pi\over 2}\,.\eeq
As a consequence, this implies that all terms in  (\ref{fpart}) involving 
$\sin\Theta$ or $\cos\Theta$ remain uniformly positive.

The lower bound $\Theta\geq \theta_0$ is an immediate consequence of (\ref{Theta2}).
To obtain an upper bound on $\Theta$, we set
$$q^\sharp~\doteq~{q(y)\over I(y)}\,.$$
By (\ref{BVP}), a differentiation yields
$$ \dot q^\sharp~=~ {c\alpha z^{\alpha-1}- q^\sharp I' \sin(\Theta) \over I}\,.$$
Next, we observe that, by (\ref{BVP}), one has
$${dz\over dq^\sharp}~=~ \ln q^\sharp\cdot 
 \cos(\Theta-\theta_0) \cdot {I \over c\alpha z^{\alpha-1}- q^\sharp I' \sin(\Theta) }~=~ \vp_1(q^\sharp)
 \cdot \ln q^\sharp
 \cdot \alpha z^{\alpha-1}\,,
 \qquad\quad\left\{\bega{rl} z(h)&=~0,\\[3mm]
q^\sharp(h)&= ~1.\enda\right.$$
In (\ref{sveg}) we can now choose $\delta \le c \alpha M^{\alpha-1}$, where $M \ge z(0)$ is an a priori bound on the mass of the stem, derived in Section \ref{s:5}. This ensures that 
$\vp_1$ is a bounded, uniformly positive function for $y$ close enough to $h$, say
$$0< c^-\leq ~\vp_1~\leq~c^+,$$
for some constants $c^-,c^+$.
Integrating, we obtain
\bel{za} z^\alpha~=~\int_0^z \alpha \zeta^{\alpha-1}\, d\zeta~
=~-\int_{q^\sharp}^1\vp_1(s)   \ln s \,ds\, =~ -\vp_2(q^{\sharp}) \int_{q^\sharp}^1 \,\ln s\, ds =~ \vp_3(q^{\sharp}) \cdot (1-q^{\sharp})^2,\eeq
and
\bel{dqs2}{dq^\sharp\over dy} ~=~{c\alpha\over\sin\Theta} \left(
-\int_{q^\sharp}^1 \vp_1(s)\,\ln s\, ds\right) ^{\alpha-1\over\alpha}~=~
\vp_4(q^\sharp) \cdot \left(
-\int_{q^\sharp}^1 \,\ln s\, ds\right) ^{\alpha-1\over\alpha}
~=~\vp_5(q^\sharp)\cdot (1-q^\sharp)^{2(\alpha-1)\over\alpha}.\eeq
Here $\vp_k$ are uniformly positive, bounded functions.
Integrating (\ref{dqs2}) we obtain
\bel{qsb}
\int_{q^\sharp}^1 {1\over\vp_5(s)} (1-s)^{2(1-\alpha)\over\alpha} \, ds~=~h-y.\eeq
To fix the ideas, assume 
$$0~<~c_3~\leq~\vp_5(s)~\le~ C_3\,.$$   Then
$$
{1\over c_3} \int_{q^\sharp}^1  (1-s)^{2(1-\alpha)\over\alpha} \, ds~=~{\alpha \over (2-\alpha) c_3} (1-q^\sharp)^{2 - \alpha \over\alpha} \, ds~\geq ~h-y.$$
\bel{qsb2}1-q^\sharp (y)~\geq~\left( (2-\alpha)c_3 \over \alpha \right)^{\alpha \over 2-\alpha } (h-y)^{\alpha \over 2-\alpha }.
\eeq
A similar argument yields
\bel{qsb3}1-q^\sharp (y)~\leq~\left( (2-\alpha)C_3 \over \alpha \right)^{\alpha \over 2-\alpha } (h-y)^{\alpha \over 2-\alpha }.
\eeq

Using $I'(y)\leq\delta$ and (\ref{qsb3}) in the equation (\ref{dpdy}) we obtain a bound of the form
\bel{py7}-p'(y)~\leq~C_1 (1-q(y))~\leq~C_2(h-y)^{\alpha \over 2-\alpha},\eeq
which yields
\bel{py4}p(y)~\leq~ {C_2\over \alpha+1}(h-y)^{2 \over 2-\alpha }.\eeq
Since $\alpha<1$, using (\ref{qsb2}) and (\ref{py4}) in (\ref{Theta2})
we obtain the limit $\Theta(y)\to \theta_0$ as $y\to h$.

On the other hand, when $y$ is bounded away from $h$, the denominator in (\ref{w})
is strictly positive and the quantity $w=w(I,p,q)$ remains uniformly bounded.
By (\ref{Theta2}), we obtain the upper bound $\Theta\leq\theta^+$, for some $\theta^+<\pi/2$.
\v
{\bf 8.} Relying on (\ref{fpart}), in this step we prove that all terms on the right hand sides of the ODEs 
in (\ref{dh1}) are uniformly integrable. 
\begi
\item[(i)] We first consider the terms appearing in the ODE for $P(y)$.
Concerning  $f_{1,p}$, as $y\to h-$ one has 
 \bel{f1p}f_{1,p}~=~\O(1)\cdot \Big(1-{q\over I}\Big)^{-1}~=~ 
 \O(1) \cdot (h-y)^{-\alpha \over 2-\alpha},\eeq
because of (\ref{qsb2}).
Since $\alpha<1$, this implies that  $f_{1,p}$ is an  integrable function of $y$. 
\item[(ii)] By the second equation in (\ref{fpart}), as $y\to h-$ one has 
\bel{f1q}
	f_{1,q}~=~\O(1)\cdot{ (1-q^\sharp) \ln(q^\sharp) \over 1-q^\sharp+q^\sharp\ln(q^\sharp) } ~= ~\O(1).
	\eeq
\item[(iii)] The term $f_2$ blows up as $y\to h-$, due to the factor $z^{\alpha-1}$.
	However, this factor is integrable in $y$ because, by \eqref{za}, \eqref{qsb2} and \eqref{qsb3}
	\begin{equation}
	z^\alpha\bigl(I(y),p(y),q(y)\bigr)~ =~ \O(1) \cdot (h-y)^{\frac{2\alpha}{2-\alpha}}.
	\label{za2}
	\end{equation}
	This implies
\bel{zeq}
f_2\bigl(I(y), p(y), q(y)\bigr)~=~\O(1)\cdot	z^{\alpha-1}\bigl(I(y),p(y),q(y)\bigr) ~= ~\O(1) \cdot (h-y)^{-1+\frac{\alpha}{2-\alpha}},
\eeq
showing that $f_2$ is integrable, because $\alpha > 0$.
\v
\item[(iv)] We now solve the linear ODE for $P$ in (\ref{dh1}) with terminal condition
$P(h)=0$. By  the estimates (\ref{f1p})-(\ref{f1q}) and (\ref{zeq}) 
one obtains a bound of the form
\bel{Phy}
P(y)~=~\O(1)\cdot (h-y)^{\alpha\over 2-\alpha},\eeq
valid in a left neighborhood of $y=h$.

\item[(v)] In a neighborhood of the origin,
the function  $f_{1,q}$ contains  a logarithm which blows up as $y \to 0+$.
However, 	this is integrable because, for $y\approx 0$,  we have
	\begin{equation*}
	\frac{q(y)}{I(y)} ~\approx ~
	\left(\diff{}{y} \frac{q(y)}{I(y)}\right) \bigg|_{y=0}\cdot  y ~=~ 
	\frac{c \alpha}{(z(0))^{1-\alpha} I(0) \sin{(\Theta(0))}} \,y,
	\end{equation*}
	and $\ln y$ is integrable in $y$.
Recalling (\ref{Iprime}), as $y>0$ ranges in a  neighborhood of the origin, we conclude
\bel{IP2}\left\{\bega{rl} I'(y) \cdot f_{1,q}f_2&=~\O(1)\cdot I'(y) f_{1,q}~=~\O(1)\cdot  y^{-\beta} \ln y,
\\[4mm]
I'(y) \cdot f_{1,p}&=~\O(1)\cdot I'(y) ~=~\O(1)\cdot  y^{-\beta} \,.\enda\right.
\eeq
This shows that, in  (\ref{dh1}), the coefficients in first equation  are uniformly integrable
in a neighborhood of the origin.
\v
\item[(vi)] 
It remains to consider the terms appearing in the ODE for $\Tilde{Q}(y)$.
We first observe that
	\begin{equation*}
	\frac{f_{2,p}}{f_2} ~= ~-\frac{\sin{\Theta}}{c \alpha} \left[1 + \frac{\alpha}{\sin^2\Theta} - 2\alpha \right] z^{-\alpha}\bigl(I(y),p(y),q(y)\bigr).
	\end{equation*}
	 As $y\to h-$, 
by (\ref{za2}) and (\ref{Phy})  this implies
\bel{f7}\frac{f_{2,p}}{f_2}\cdot P~= ~\O(1)\cdot (h-y)^{-2\alpha\over 2-\alpha}\cdot (h-y)^ {\alpha\over 2-\alpha},
\eeq
which is integrable for $\alpha < 1$.  
\v
\item[(vii)]
	Finally, as $y\to h-$, we consider
	\bel{f8}\bega{rl}\ds
	\frac{f_{2,I}}{f_2} &=~\ds -\frac{\sin{\Theta}}{c \alpha}\left[\frac{(1-\alpha) \sin{\theta_0}}{\sin^2\Theta} + \frac{\sin\left(\Theta-\theta_0\right)}{\cos\Theta}\left(1 + \frac{\alpha}{\sin^2\Theta} -2\alpha \right)  \right] \frac{1-\frac{q}{I(y)}}{z^{\alpha}\bigl(I(y),p(y),q(y)\bigr)}
	\\[4mm]
	&=~\ds \O(1)\cdot (1-q^\sharp) z^{-\alpha}\bigl(I(y),p(y),q(y)\bigr)
	~=~\O(1) \cdot (h-y)^{\alpha\over 2-\alpha} \cdot(h-y)^{-2\alpha\over 2-\alpha} ,
	\enda
	\eeq
	which  is integrable in $y$ since $\alpha < 1$.
Similarly,  by (\ref{f7}), (\ref{dpdy}), and (\ref{qsb3}), it follows
\bel{f10}{f_{2,p} \over f_2}\cdot f_1~=~   \O(1)\cdot (h-y)^{-2\alpha\over 2-\alpha} \cdot (h-y)^{\alpha\over 2-\alpha},
\eeq
which is again integrable.
\endi	
\v
{\bf 9.}  The proof can now be accomplished by a contradiction argument.
If the conclusion of the theorem were not true, one could find a sequence
of absolutely continuous, non-decreasing functions $I_n:\R_+\mapsto [0,1]$,
all satisfying (\ref{Iprime}), with $I_n(0)\to 1$, and such that, for each $n\geq 1$,
the optimization problem {\bf (OP2)} has two distinct solutions, say
$(\check \theta_n, \check u_n)$ and $(\hat\theta_n, \hat u_n)$.
As a consequence, for each $n\geq 1$ the system  (\ref{ODE20}) has two solutions.
To fix the ideas, let the first solution be defined on $[0, \check h_n]$ and the second on $[0,
\hat h_n]$, with $\check h_n<\hat h_n$. 
These two solutions will be denoted by $(\check p_n, \check q_n, \check z_n)$ and $(\hat p_n, \hat q_n, \hat z_n)$.
They both satisfy the boundary conditions
\bel{bcc}\check {p}_n(\check {h}_n)\, =\, \hat{p}_n(\hat{h}_n) \,=\, 0,
\qquad \check q_n(\check h_n) \,=\, I(\check h_n),\qquad \hat q_n(\hat h_n) \,=\, I(\hat h_n),
\qquad \check {q}_n(0)\, =\, \hat{q}_n(0) \,=\,0.\eeq
Thanks to the last identity, by the mean value theorem
there exists some intermediate point  $k_n\in [\check {h}_n,\hat{h}_n]$ such that, with the notation introduced at (\ref{pqyh}),
 \bel{parqn}\pdiff{q_n}{h}(0,k_n) ~=~ 0.\eeq
For each $n\geq 1$ consider the corresponding system
 \bel{dhn}\left\{\bega{rl}
	P_n'(y) &=\,\ds- I'_n(y) \left[ f_{1,p} P_n + f_{1,q} f_2 \Tilde{Q}_n \right], \\[3mm]
	\Tilde{Q}'(y) &=~\ds \frac{f_{2,p}}{f_2} P_n
	 - \frac{I_n'(y) [f_{2,I} - f_{2,p} f_1]}{f_2} \Tilde{Q}_n,
	\enda\right.\qquad\qquad \left\{ \bega{rl} P_n(k_n)&=~0,\\[3mm]
	\Tilde Q_n(k_n)&=~-1.\enda\right.\eeq
Since $f_2\bigl(I_n(0),p_n(0,k_n),0\bigr) > 0$, by (\ref{parqn}) it follows  
\bel{TQ0}\Tilde Q_n (0)~=~{1\over f_2\bigl(I_n(0),p_n(0, k_n),0\bigr)}\cdot \pdiff{q_n}{h}(0,k_n) ~=~0.\eeq
Let
	\begin{equation*}
	P_n(y)~ \doteq ~\pdiff{p(y,k_n)}{h},\qquad\qquad \Tilde{Q}_n(y) ~\doteq 
	~{1\over f_2\bigl(I_n(y),p_n(y, k_n),q_n(y, k_n)\bigr)}\cdot\pdiff{q(y,k_n)}{h},
	\end{equation*}
	be the solutions to (\ref{dhn}).  By the previous steps, their derivatives $\bigl(P_n',  \Tilde{Q}_n' \bigr)_{n \ge 1}$
form a sequence of uniformly integrable functions  defined on the   intervals $[0,k_n]$. 
Note that the existence of an upper bound $\sup_n k_n\,\doteq\,h^+<+\infty$  follows  from the existence proof.

Thanks to the uniform integrability, by possibly taking a subsequence, 
we can assume the convergence $k_n\to \ov h\in [0,h^+]$,  the weak convergence of derivatives 
$P'_n\wto P'$, $\Tilde Q_n' \rightharpoonup  \Tilde{Q}'$ in $\L^1$, and the convergence 
$$P_n ~\to ~P,\qquad \Tilde Q_n~\to~ \Tilde{Q},$$
uniformly on every subinterval $[0,h]$ with $h<\bar h$.

Recalling that every $I_n'$ satisfies the uniform bounds  (\ref{Iprime}), since $I_n(y)\to I(y)\equiv 1$ uniformly
for all $y\geq 0$, we conclude that $(P,\Tilde Q)$  provides a solution to the linear system (\ref{dh1}) 
on $[0, \bar h]$, corresponding to the constant function $I(y)\equiv 1$.
 We now observe that, when $I(y)\equiv 1$, the solution to (\ref{dh1}) is  $P(y)\equiv 0$ and $\Tilde Q(y)\equiv -1$.
On the other hand, our construction yields
	\begin{equation*}\Tilde Q(0)~=~ \lim_{n\to \infty}\Tilde{Q}_n(0) ~=~0.
	\end{equation*}
	This contradiction achieves the proof of Theorem~\ref{t:61}.
	\endproof

\section{Existence of an equilibrium solution}
\label{s:6}
\setcounter{equation}{0}
Given a nondecreasing light intensity function $I:\R_+\mapsto [0,1]$, in the previous section
we proved the existence of an optimal  solution 
$(\theta^*, u^*)$ for the maximization problem {\bf (OP2)}.    

Conversely, let $\rho_0>0$ be the constant density of stems, i.e.~the number of stems
growing per unit area. 
If all stems have the same configuration, described by the couple of functions $y\mapsto (\theta(y), u(y))$ as in
(\ref{max33}), then the corresponding intensity of light at height $y$ above ground is computed as
\bel{Iy}I^{(\theta,u)}(y)~\doteq~\exp\left\{ -{\rho_0\over\cos\theta_0} \int_y^{+\infty} {u(\zeta)\over \sin\theta(\zeta)}\, d\zeta\right\}.
\eeq

The main goal of this section is to find a competitive equilibrium, i.e.~a fixed point of the composition of the two maps
$I\mapsto (\theta^*,u^*)$ and $(\theta,u)\mapsto I^{(\theta,u)}$.

\v
\begin{definition}\label{d:E2} Given an angle $\theta_0\in \,]0, \pi/2[\,$ and a constant $\rho_0>0$, we say 
that 
the light intensity function $I^*:\R_+\mapsto [0,1]$ and the stem configuration 
$(\theta^*, u^*) :\R_+\mapsto 
[\theta_0, \pi/2]\times \R_+$
yield a {\bf competitive equilibrium} if the following holds.
\begi
\item[(i)] The couple $(\theta^*, u^*)$  provides an optimal solution to the optimization problem 
{\bf (OP2)}, with light intensity function $I=I^*$.
\item[(ii)] The identity  $I^*= I^{(\theta^*, u^*)}$ holds.\endi
\end{definition}

The main result of this section provides the existence of a competitive equilibrium, 
assuming that the density $\rho_0$ of stems is sufficiently small.

\begin{theorem}\label{t:6} Let an angle $\theta_0\in \,]0, \pi/2[\,$ be given. Then,  for all $\rho_0>0$
sufficiently small, a unique  competitive equilibrium $(I^*, \theta^*, u^*)$ exists.
\end{theorem}

{\bf Proof.} {\bf 1.} Setting $C=1$ and $\beta= 1/2$ in (\ref{Iprime}), we define the family of functions 
\bel{Fdef}
\F~\doteq~\big\{ I:\R_+\mapsto [1-\delta,\,1]\,;\quad I ~\hbox{ is absolutely continuous},\quad
 I'(y)\in \bigl[0,\,  y^{-1/2} \,\bigr]
\quad\hbox{for a.e.~} y>0\Big\},\eeq
where  $\delta>0$ is chosen small enough so that the conclusion of Theorem~\ref{t:61}
holds.
\v
{\bf 2.}
For each $I\in \F$, let $(\theta^{(I)},u^{(I)})$ describe the corresponding 
optimal stem.  
Calling 
$$h^{(I)}~=~\sup~\bigl\{ y\geq 0\,;~~u^{(I)}(y)>0\bigr\}$$
the height of this stem, by the a priori bounds proved in Section~\ref{s:55} we have a uniform bound
$$h^{(I)}~\leq~h^+$$
for all $I\in \F$.   
Let $p^{(I)}, q^{(I)}:[0, h^{(I)}]\mapsto \R_+$ be the corresponding solutions of (\ref{bvp2}).
For convenience, we extend all these functions to the  larger interval $[0, h^+]$ by setting
$$p^{(I)}(y)\,\doteq\, p^{(I)}\bigl(h^{(I)}\bigr), \qquad\qquad  q^{(I)}(y)\,\doteq\, q^{(I)}\bigl(h^{(I)}\bigr)
\qquad \forall y\in [h^{(I)}, h^+].
$$
\v
{\bf 3.}
By the estimates proved in Section~\ref{s:55}, if we choose $\rho_0>0$ small enough, it follows that
the corresponding light intensity function $I^{(\theta,u)}$  at (\ref{Iy}) is again in $\F$.
A competitive equilibrium will be obtained by 
constructing a fixed point of the composition of the two maps
\bel{IUI}\Lambda_1:
I~\mapsto~\bigl(\theta^{(I)},\,u^{(I)}\bigr),\qquad\qquad \Lambda_2:(\theta,u)~\mapsto~I^{(\theta,u)}.\eeq
In order to use Schauder's theorem, we need to check the continuity of these map, in a suitable topology.

We start by observing that 
$\F\subset\C^0([0,h^+])$ is a compact, convex set.  
Again by the analysis in Section~\ref{s:55}, as  $I$ varies within the domain $\F$, 
the corresponding functions
$\theta^{(I)}$ are uniformly bounded in $\L^\infty([0,h^+])$, while $u^{(I)}$
is uniformly bounded in $\L^1  ([0,h^+])$.

{} From the estimate (\ref{py7}) it follows that the functions $p^{(I)}$ are equicontinuous on $[0, h^+]$.
Recalling that $q= q^\sharp\cdot I$,  by (\ref{dqs2}) we conclude that the functions $q^{(I)}$
are equicontinuous as well.
\v
{\bf 4.} By the analysis in Section~\ref{s:55}, for any $I\in \F$,  the solution to the 
system of optimality conditions (\ref{BVP}) satisfies
\bel{TUP2}
\theta_0~\leq~\Theta(I(y), p(y), q(y))\,\leq~\theta^+\,,\qquad\qquad  
 c_0\,y\,\leq\, {q(y)\over I(y)}~\leq~1,\eeq
 for some $\theta^+<\pi/2$ and $c_0>0$ sufficiently small. In view of (\ref{U}), this implies
 \bel{TUP3}U(I(y), p(y), q(y))~\doteq~-\ln\left(q(I)\over I(y)\right)\,\cos\bigl(\Theta(I(y), p(y), q(y))-\theta_0
 \bigr)~\leq~-\ln (c_0 y).
\eeq
Motivated by (\ref{TUP2})-(\ref{TUP3}), we consider the set of functions
\bel{Udef}\U~\doteq~\Big\{ (\theta, u)\in \L^1\bigl([0, h^+]\,;~\R^2\bigr),\quad \theta(y)\in [\theta_0, \theta^+],
~~0\leq u(y)\leq -\ln (c_0 y)\Big\}.\eeq

Thanks to the uniform bounds imposed on $\theta$ and $u$ in the definition (\ref{Udef}), 
the  continuity of the map $\Lambda_2:\U\mapsto \C^0$, defined at (\ref{Iy})
is now straightforward.  

\v
{\bf 5.} To prove the continuity of the map $\Lambda_1$,  consider a sequence of functions $I_n\in \F$,
with $I_n\to I$ uniformly on $[0, h^+]$.    Let $(\theta_n, u_n): [0, h^+]\mapsto \R^2$
be the corresponding unique optimal solutions.

We claim that $(\theta_n, u_n)\to (\theta,u)$ in $\L^1([0, h^+])$, where 
$(\theta, u)$ is the unique optimal  solution, given the light intensity $I$. 

To prove the claim, let $(p_n, q_n)$ be the corresponding solutions of the system (\ref{bvp2}).
By the estimates on $p', q'$ proved in Section~\ref{s:55}, the functions $(p_n,q_n)$ are equicontinuous.
From any subsequence we can thus extract a further subsequence and obtain the convergence
\bel{pqI}p_{n_j}\to \Hat p,\qquad q_{n_j}\to \Hat q,\qquad I_{n_j}\to I,\eeq
for some functions $\Hat p,\Hat q$,
uniformly on $[0, h^+]$.

For every $j\geq 1$ we now have
$$\theta_{n_j}(y)~=~\Theta\bigl(I_n(y), p_n(y), q_n(y)\bigr), \qquad\qquad
u_{n_j}(y)~=~U\bigl(I_n(y), p_n(y), q_n(y)\bigr), $$
where $U$ and  $\Theta$ are the functions in (\ref{U})-(\ref{Theta2}).
By the dominated convergence theorem, the  convergence (\ref{pqI}) together with the uniform integrability of $\theta_{n_j}$
 and $u_{n_j}$
yields the $\L^1$ convergence
\bel{L1c} \|\theta_{n_j}-\Hat \theta\|_{\L^1}~\to ~0,\qquad\qquad \|u_{n_j}-\Hat u\|_{\L^1}~\to ~0.\eeq
In turn this implies that $(\Hat p,\Hat q)$ provide a solution to the 
problem (\ref{bvp2}), in connection with the light intensity $I$.  
 By uniqueness, $\Hat p=p$ and $\Hat q=q$.  Therefore,
$\Hat \theta=\theta$ and $\Hat u=u$ as well.

The above argument shows that, from any subsequence, one can extract a further subsequence
so that the $\L^1$-convergence (\ref{L1c}) holds.  Therefore, the entire sequence $(\theta_n,u_n)_{n\geq 1}$
converges to $(\theta,u)$ in $\L^1([0,h^+])$.   This establishes the continuity of the map $\Lambda_1$.
\v
{\bf 6.} The map $\Lambda_2\circ\Lambda_1$ is now a continuous map of the 
compact, convex domain $\F\subset\C^0([0, h^+])$ into itself.   By Schauder's theorem it admits
a fixed point $I^*(\cdot)$.   By construction, the optimal stem configuration 
$\bigl(\theta^{(I^*)}, u^{(I^*)}\bigr)$ yields a competitive equilibrium, in the sense of Definition~\ref{d:E2}.
\v
{\bf 7.} To prove uniqueness, we derive a set of necessary conditions satisfied by the 
equilibrium solution, and show that this system has a unique solution.

Using (\ref{U}) and (\ref{trig1}), we can rewrite the light intensity function (\ref{Iy}) as
\begin{equation*}
	I(y)~=~\exp\bigg\{\frac{\rho_{0}}{\cos\theta_{0}}\int_{y}^{\infty}\ln\Big(\frac{q}{I}\Big)\frac{1+w\sin\theta_{0}}{\sin\theta_{0}+w}\,d\zeta\bigg\},
\end{equation*}
where $w=w(I,p,q)$ is the function introduced at (\ref{w}).
Differentiating w.r.t.~$y$ one obtains
\bel{Iprim}
	I'(y)~=\,-\frac{\rho_{0}}{\cos\theta_{0}}\ln\bigg(\frac{q}{I}\bigg)\frac{1+w\sin\theta_{0}}{\sin\theta_{0}+w}\cdot I~\doteq~ f_{3}(I,p,q).
\eeq
Combining (\ref{Iprim}) with (\ref{bvp2}), we conclude that the competitive equilibrium  
satisfies the system of equations and boundary conditions
\bel{ceq}
	\left\{\bega{rl}\ds p'(y)&=- f_1\bigl(I(y),p(y),q(y)\bigr)\cdot f_{3}(I(y),p(y),q(y)), \\[3mm]
	q'(y) & =~f_2\bigl(I(y),p(y),q(y)\bigr),\\[3mm]
	I'(y)&=~f_{3}(I(y),p(y),q(y)),
	\enda \right.\qquad\qquad \left\{
	\bega{rl}
	p(h) &= ~0,\\[2mm]
	q(h) &= ~1,\\[2mm] I(h)&=~1,\enda \right.\qquad\quad
\eeq
together with 
\bel{q00}
	q(0)~=~0.\eeq
Here the common height of the stems $h>0$ is a constant to be determined.
\v
{\bf 8.} The uniqueness of solutions to (\ref{ceq}) will be 
achieved by a contradiction argument.
Since this is very similar to the one used in the proof of Theorem~\ref{t:61}, we only 
sketch the main steps.

In analogy with (\ref{PQdef}), (\ref{TQdef}), denote by $p(y,h), q(y,h), I(y,h)$ the unique solution
to the Cauchy problem (\ref{ceq}), with terminal conditions given at $y=h$.  Consider the functions
\begin{equation*}
	P(y)\,\doteq\,\frac{\partial p(y,h)}{\partial h}, \qquad \Tilde{Q}(y)\,\doteq\,\frac{1}{f_{2}(I,p,q)}\frac{\partial q(y,h)}{\partial h}, \qquad J(y)\,\doteq\,\frac{\partial I(y,h)}{\partial h}\,.
\end{equation*}
By (\ref{ceq}), these functions satisfy
\begin{equation}
\label{syst2}
\left\{\bega{rl}
	P'(y)&=\,-\big[f_{3,I}f_{1}+f_{3}f_{1,I}\big]J-\big[f_{3,p}f_{1}+f_{3}f_{1,p}\big]P-\big[f_{3,q}f_{1}+f_{3}f_{1,q}\big]f_{2}\Tilde{Q},\\[4mm]
	\Tilde{Q}'(y)&=~\ds\frac{f_{2,I}}{f_{2}}J+\frac{f_{2,p}}{f_{2}}P-\frac{f_{3}}{f_{2}}\big[f_{2,I}-f_{2,p}f_{1}\big]\Tilde{Q},\\[4mm]
	J'(y)&=~\ds f_{3,I}J+f_{3,p}P+f_{3,q}f_{2}\Tilde{Q},
\enda\right.
\end{equation}
with boundary conditions
\begin{equation*}
	P(h)=0, \qquad \Tilde{Q}(h)=-1, \qquad J(h)=0.
\end{equation*}
Set $d_{0}=\frac{\rho_{0}}{\cos\theta_{0}}$. Several of the partial derivatives 
on the right-hand side of \eqref{syst2} were computed in (\ref{fpart}). 
The remaining ones are
\begin{align*}
	f_{1,I}(I,p,q)&=\frac{q}{I^{2}}\cdot\frac{1+w\sin\theta_{0}}{\sin\theta_{0}+w}-\frac{\cos^{2}\theta_{0}}{(\sin\theta_{0}+w)^{2}}\frac{w^{2}}{p}\Big[1-\frac{q}{I}\Big],\\
	f_{3,I}(I,p,q)&=-d_{0}\bigg[\Big(\ln\Big(\frac{q}{I}\Big)-1\Big)\frac{1+w\sin\theta_{0}}{\sin\theta_{0}+w}-I\ln\Big(\frac{q}{I}\Big)\frac{\cos^{2}\theta_{0}}{(\sin\theta_{0}+w)^{2}}\frac{w^{2}}{p}\Big(1-\frac{q}{I}\Big)\bigg],\\
	f_{3,p}(I,p,q)&=d_{0}I\ln\Big(\frac{q}{I}\Big)\frac{\cos^{2}\theta_{0}}{(\sin\theta_{0}+w)^{2}}\frac{w}{p},\\
	f_{3,q}(I,p,q)&=-d_{0}I\bigg[\frac{1}{q}\cdot\frac{1+w\sin\theta_{0}}{\sin\theta_{0}+w}+\Big[\ln\Big(\frac{q}{I}\Big)\Big]^{2}\frac{\cos^{2}\theta_{0}}{(\sin\theta_{0}+w)^{2}}\frac{w^{2}}{p}\bigg].
\end{align*}
By the same arguments used in step {\bf 8}  of the proof of Theorem~\ref{t:61}, we
conclude that the right-hand side of \eqref{syst2} is uniformly integrable.
\v
{\bf 9.} Let  a density $\rho_0>0$ be given.
Assume that the problem (\ref{ceq})-(\ref{q00}) has two distinct solutions
% $(\hat{p}_{\rho_{0}},\hat{q}_{\rho_{0}},\hat{I}_{\rho_{0}})$ and $(\check{p}_{\rho_{0}},\check{q}_{\rho_{0}},\check{I}_{\rho_{0}})$, defined on $[0,\hat{h}_{\rho_{0}}]$ and $[0,\check{h}_{\rho_{0}}]$  say with $\hat h_{\rho_0}<\check h_{\rho_0}$.
$(\hat{p},\hat{q},\hat{I})$ and $(\check{p},\check{q},\check{I})$, defined on $[0,\hat{h}]$ and $[0,\check{h}]$  say with $\hat h<\check h$.
 Since $\hat q(0)=\check q(0)=0$, by the mean value theorem
 there exists  $k\in [\hat h,\check h]$ such that ${\partial q\over\partial h}(0,k)=0$.

If multiple solutions exist for arbitrarily small values of the density $\rho_0$,
we can find a decreasing sequence $\rho_{0,n}\downarrow 0$ and corresponding
solutions $P_n, Q_n, I_n$ of (\ref{syst2}), defined for $y\in [0, k_n]$, such that
\bel{bcn}
	P_n(k_n)=0, \qquad \Tilde{Q}_n(k_n)=-1, \qquad J_n(k_n)=0,
	\qquad \Tilde{Q}_n(0)=0.
\eeq
Thanks to  the uniform integrability of the right hand sides of (\ref{syst2}), by possibly extracting a subsequence we can achieve the  convergence
 $k_n\rightarrow \bar{h}\in[0,h^{+}]$, the weak convergence $P_n'\rightharpoonup P'$, $\Tilde{Q}_n'\rightharpoonup \Tilde{Q}'$, $J_n'\rightharpoonup J'$ in $\L^1$,
and the strong convergence
\begin{equation*}
	P_n\rightarrow P, \qquad \Tilde{Q}_n\rightarrow \Tilde{Q}, \qquad J_n\rightarrow J,
\end{equation*}
uniformly on every subinterval $[0,h]$ with $h<\bar{h}$. 

To reach a contradiction, we observe that
\begin{equation*}
	J_n(y)~=~-\int_{y}^{k_n}J_n'(z)\,dz
\end{equation*}
and the right-hand side of $J_n'$ in \eqref{syst2} consists of uniformly integrable terms which are multiplied by $\rho_{0,n}$.
This implies  $J(y)\equiv 0$. This corresponds to the case of an intensity function $I(y)\equiv 1$. But in this case we know that $\Tilde{Q}(y)\equiv-1$,  contradicting the fact that, by (\ref{bcn}),
\begin{equation*}
	\Tilde{Q}(0)~=~\lim_{n\to\infty}\Tilde{Q}_n(0)~=~0.
\end{equation*}
\endproof
\section{Stem competition on a domain with boundary}
\label{s:7}
\begin{figure}[ht]
\centerline{\hbox{\includegraphics[width=6cm]{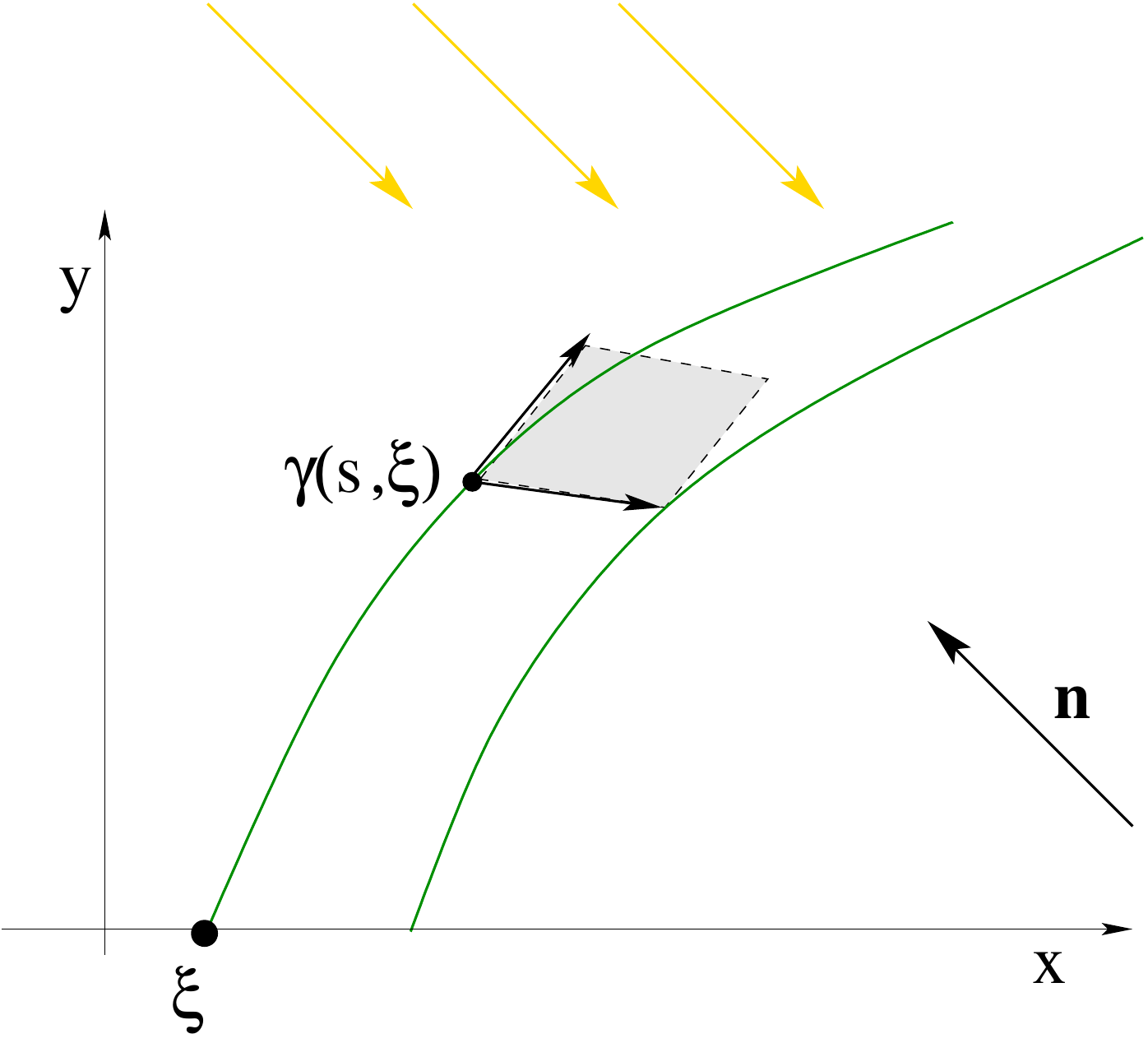}\qquad
\includegraphics[width=8.5cm]{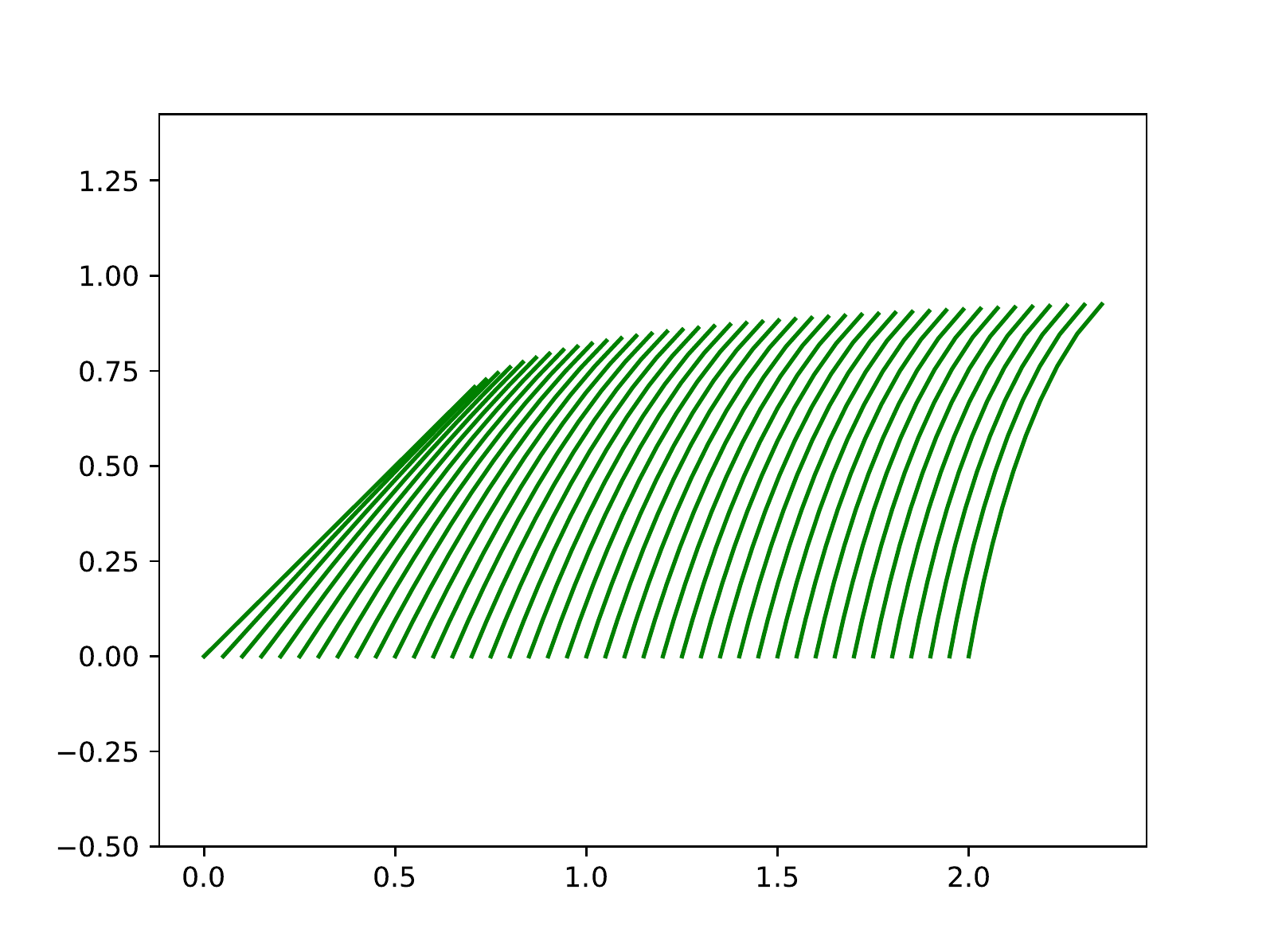}}      }
\caption{\small  Left: to leading order, the amount of vegetation in the shaded region is proportional to 
$\kappa\,\bar \rho(\xi)d\xi d s$.   Since the area is computed in terms of the cross product
${\partial\gamma\over\partial\xi}\times {\partial\gamma\over\partial s}$,
this motivates the formula  (\ref{vegd}).  Right: a possible competitive equilibrium, 
where the light rays come from the direction  $\bfn=({-1\over\sqrt 2}, {1\over\sqrt 2})$ and stems are
distributed along the positive half line, with density as in (\ref{hls}).  
 In this case, stems originating from points close to the origin have no incentive to grow upward, 
 because they already receive a nearly maximum light intensity.   Hence they 
 bend to the right, almost perpendicularly to the light rays.    }
\label{f:pg44}
\end{figure}

\setcounter{equation}{0}
We consider here the same model introduced in Section~2, 
where all stems have fixed length $\ell$ and constant thickness $\kappa$.  But we now  allow the
sunlight intensity $I=I(x,y)$ to vary w.r.t.~both  variables $x,y$. As  shown
in Fig.~\ref{f:pg44}, left,
we denote by
 \bel{stemp}s~\mapsto ~\gamma(s,\xi)~=~(x(s), y(s)),\qquad s\in [0,\ell],\eeq
the arc-length parameterization of the stem whose root is located at $(\xi,0)$, and write $g$ for the function introduced at (\ref{g}).
This leads to the optimization problem
 \begi
\item[{\bf (OP3)}] {\it Given a  light intensity function $I=I(x,y)$, 
find a control
$s\mapsto\theta(s)\in [0, \,\pi]$ which maximizes the integral
\bel{I6} \int_0^\ell I(x(s),y(s))\, g(\theta(s))\, ds\eeq
subject to}
\bel{ic3} {d\over ds} (x(s), y(s))~=~(\cos\theta(s), \sin\theta(s)), 
\qquad\qquad (x(0), y(0))\,=\,(\xi,0).\eeq
\endi

Next, consider a function $\bar\rho(\xi)\geq 0$ describing the density of stems 
which grow near  $\xi\in\R$. 
At any point in space reached by a stem, i.e. such that
$$(x,y)=\gamma(s,\xi)\qquad\hbox{for some}~~\xi\in\R,~~s\in [0,\ell],$$
the density of vegetation is %(see Fig.~\ref{f:pg44})
\bel{vegd}
\rho(x,y)~=~\rho(\gamma(s,\xi))~=~\kappa\,
\bar\rho(\xi)\cdot\left[{\partial\gamma\over\partial\xi}\times {\partial\gamma\over\partial s}
\right]^{-1}.\eeq
The light intensity at a point $P=(x,y)\in \R^2$ is now given by
\bel{IP}
I(P)~=~\exp\left\{-\int_0^{+\infty} \rho(P+t\bfn)\, dt\right\}.\eeq
\v
\begin{definition}
Given the constants $\ell,\kappa$ and the density $\bar\rho\in \L^\infty(\R)$,  we say that the maps
$\gamma:[0,\ell]\times\R$ and $I:\R\mapsto\R_+\mapsto [0,1]$ yield 
a {\bf competitive equilibrium} if the following holds.
\begi
\item[(i)]  For each $\xi\in\R$, the stem $\gamma(\cdot,\xi)$ provides an optimal solution to {\bf (OP3)}.
\item[(ii)] The function $I(\cdot)$ coincides with the light intensity determined by (\ref{vegd})-(\ref{IP}).
\endi\end{definition}
\v
We shall not analyze the existence or uniqueness of the competitive equilibrium, in the case 
where the distribution of stem roots is not uniform. 
We only observe that, if the stem $\gamma(\cdot,\xi)$ in (\ref{stemp}) is optimal, 
the necessary conditions yield the existence of a dual vector
$s\mapsto \bfp(s)$ satisfying
\bel{nc9}\dot \bfp(s)~=~- \nabla I\bigl(x(s),y(s)\bigr)\, g(\theta(s)),\qquad\qquad \bfp(\ell)\,=\,(0,0),\eeq
and such that, for a.e.~$s\in [0,\ell]$, the optimal angle $\theta^*(s)$ satisfies
\bel{oc4}\theta^*(s)~=~\argmax_{\theta} \Big\{ \bfp(s)\cdot(\cos \theta, \sin\theta) +I(x(s), y(s)) g(\theta)\Big\}.\eeq
Differentiating the expression on the right hand side of (\ref{oc4}) %and using (\ref{nc9}) 
one obtains
an implicit equation for $\theta^*(s)$, namely
\bel{nc44}
I\bigl(x(s),y(s)\bigr))g'(\theta^*(s))+\bfp(s)
\cdot \bfn(s)~=~0\eeq
for a.e.~$s\in [0,\ell]$. Here $\bfn(s)~\doteq~\bigl( - \sin\theta(s), \cos\theta(s)\bigr)$
is the unit vector perpendicular to the stem. Moreover, by (\ref{nc9}) one has
$$\bfp(s)~=~ \int_s^\ell   \nabla I\bigl(x(\sigma), y(\sigma))\, g(\theta^*(\sigma)\bigr)\, d\sigma.$$
An interesting case is where stems grow only on the half line $\{\xi\geq 0\}$. For example, one can take
\bel{hls} \bar\rho(\xi)~=~\left\{\bega{cl}
0\quad &\hbox{if}~~\xi<0, \cr
b^{-1}\xi\quad &\hbox{if}~~\xi\in [0, b],\cr
 1\quad &\hbox{if}~~\xi> b.\enda\right.\eeq
 In this case, we conjecture that the competitive equilibrium has the form illustrated in 
 Fig.~\ref{f:pg44}, right.

\section{Concluding remarks}

A motivation for the present study was to understand whether competition
for sunlight could explain phototropism, i.e.~the tendency
of plant stems to bend toward the light source.
A naive approach may suggest that, if a stem bends in the direction of the light rays, 
the leaves will be closer to the sun and hence gather more light.  
However, since the average distance of the earth from the sun
is approximately 90 million miles, 
getting a few inches closer cannot make a difference.  

\begin{figure}[ht]
\centerline{\hbox{\includegraphics[width=8cm]
{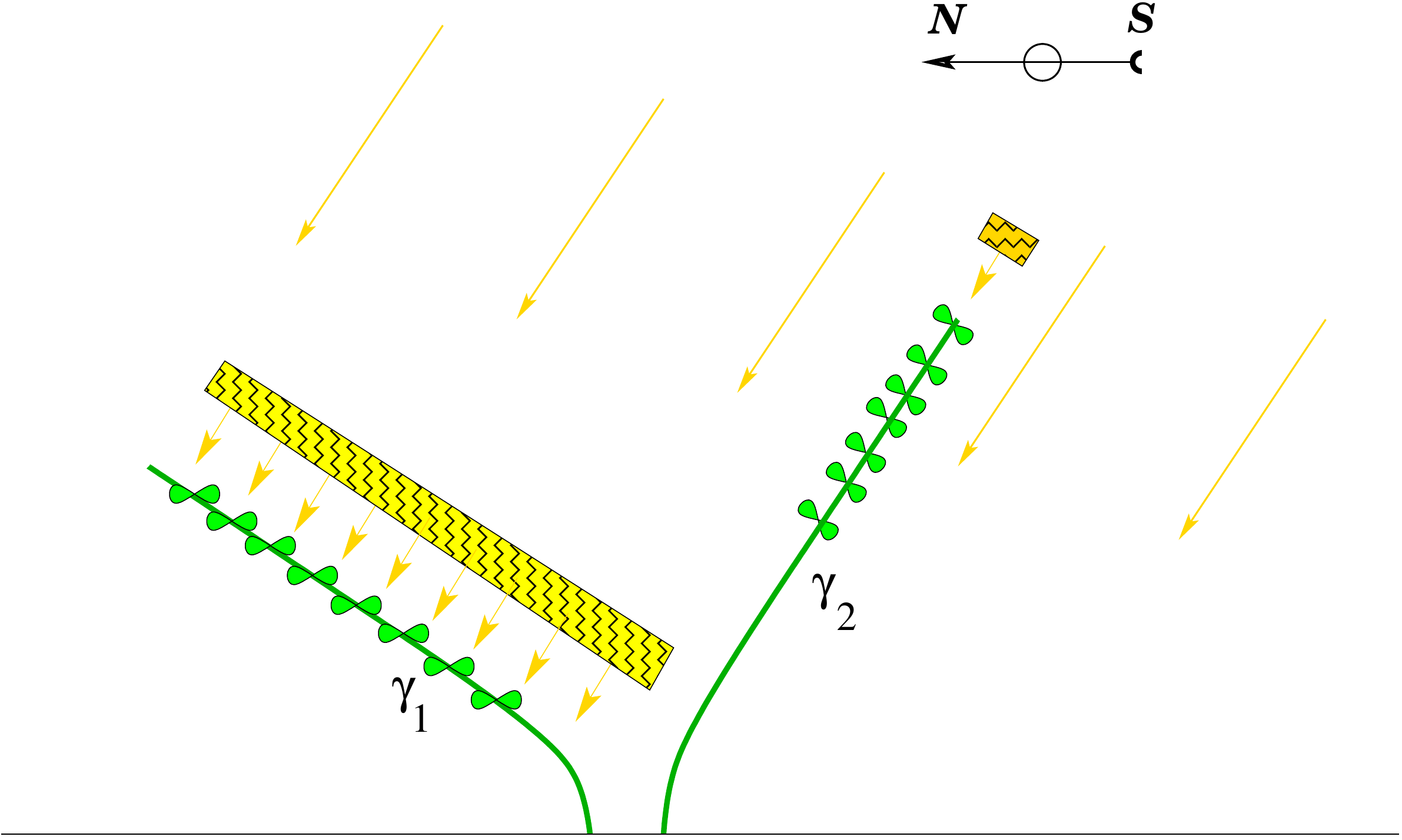}}}
\caption{\small  The stem $\gamma_1$, oriented perpendicularly to the sun rays,
collects much more sunlight than $\gamma_2$.
Indeed, $\gamma_1$ would give the best orientation for solar panels. Notice that $\gamma_2$
minimizes the sunlight gathered because the upper leaves put the lower ones in shade.}
\label{f:pg21}
\end{figure}

As shown in Fig.~\ref{f:pg21}, if a single stem were present,
to maximize the collected sunlight it should be perpendicular to the light rays, not 
parallel.  
In the presence of competition
among several plant stems, our analysis shows that  the best configuration
is no longer perpendicular to light rays:  the lower part of the stems should grow in a nearly vertical direction,
while the upper part bends away from the sun.

Still, our competition models do not predict the tilting of stems in the direction of the sun rays.  
This may be due to the fact that these models are ``static", i.e., they do not describe how plants grow
in time.    This leaves open the possibility of introducing further models that can explain  phototropism 
in a time-dependent framework.
 As suggested in~\cite{RLP}, the
preemptive conquering of space, in the direction of the light rays, can be an advantageous strategy.  
We leave these issues for future investigation.

\vs
{\bf Acknowledgment.} The research of
A.Bressan was partially supported by NSF, with grant DMS-1714237, 
``Models of controlled biological growth".
S.T.Galtung was supported in part by a grant from the U.S.-Norway Fulbright Foundation.
A.Reigstad was supported by the grant ``Waves and Nonlinear Phenomena" 
from the Research Council of Norway. S.T.Galtung and A.Reigstad are very grateful to the Department of Mathematics at Penn State University for the generous hospitality during the academic year 2018/2019.

\v

\end{document}